\documentclass[leqno,11pt]{amsart}

\usepackage{amsmath,amstext,amssymb,amsopn,amsthm,mathrsfs}
\usepackage{verbatim}
\usepackage{enumerate}

\newcommand*\mr{\mathbb{R}}

\newcommand*\Lafi{L^{\alpha}_{\varphi}}
\newcommand*\Lal{L^{\alpha}_{\ell}}

\newcommand*\Safi{\mathcal{S}_{\varphi}^{\alpha}}
\newcommand*\Sal{\mathcal{S}_{\ell}^{\alpha}}
\newcommand*\Safid{(\mathcal{S}_{\varphi}^{\alpha})'}
\newcommand*\Se{\mathcal{S}_e}
\newcommand*\Bifi{\dot{B}^{\sigma,\Lafi ,m}_{p,q}}
\newcommand*\TiLifi{\dot{F}^{\sigma,\Lafi ,m}_{p,q}}
\newcommand*\Patm{P^{\alpha}_{t,m}}
\newcommand*\Pam{\tilde{P}^{\alpha}_{t,m}}
\newcommand*\patm{p^{\alpha}_{t,m}}
\newcommand*\pat{p^{\alpha}_{t}}
\newcommand*\Pat{P^{\alpha}_{t}}

\title[Besov and Triebel-Lizorkin spaces]
{Besov and Triebel-Lizorkin spaces associated with Laguerre expansions of Hermite type}

\author[P{.} Plewa]{Pawe\l{} Plewa}
\address{Pawe\l{} Plewa \newline
			Faculty of Pure and Applied Mathematics, 
      Wroc\l{}aw University of Science and Technology       \newline
      Wyb{.} Wyspia\'nskiego 27,
      50--370 Wroc\l{}aw, Poland      
      }
\email{pawel.plewa@pwr.edu.pl}

\allowdisplaybreaks

\theoremstyle{plain}
\newtheorem{thm}{Theorem}[section]

\newtheorem{lm}[thm]{Lemma}
\newtheorem{prop}[thm]{Proposition}
\newtheorem{cor}[thm]{Corollary}
\newtheorem{defin}[thm]{Definition}

\theoremstyle{definition}

\theoremstyle{remark}
\newtheorem*{rem*}{Remark}
\newtheorem{rem}[thm]{Remark}

\theoremstyle{plain}

\begin{document}

\begin{abstract}
Homogeneous Besov and Triebel-Lizorkin spaces associated with multi-dimensional Laguerre function expansions of Hermite type with index $\alpha\in [-1/2,\infty)^d\setminus(-1/2,1/2)^d$, $d\geq 1$, are defined and investigated. To achieve expected goals Schwartz type spaces on $\mr^d_+$ are introduced and then tempered type distributions are constructed. Also, ideas from a recent paper of Bui and Duong on Besov and Triebel-Lizorkin spaces associated with Hermite functions expansions are used. This means, in particular, using molecular decomposition and an appropriate form of a Calder\'on reproducing formula.

\end{abstract}

\maketitle
\footnotetext{
The paper was author's master thesis written under the supervision of Professor Krzysztof Stempak.

}

\section{Introduction}
In this paper we explore the homogeneous Besov and Triebel-Li\-zor\-kin spaces (resp. $\Bifi$ and $\TiLifi$) in terms of the Laguerre function expansions of Hermite type for $0<p<\infty,\ 0<q\leq\infty$ (for the Besov space we allow $p=\infty$), $\alpha\in [-1/2,\infty)^d\setminus(-1/2,1/2)^d$ and $\sigma\in\mr$. Note that the classical theory of these function spaces as the spaces associated to Laplacians or their square roots on $\mr^d$ is an important part of the function spaces theory due to its applications in harmonic analysis and the theory of partial differential equations.

The range of the admissible Laguerre type multi-index $\alpha$ in great part of this paper is full, that is $\alpha\in (-1,\infty)^d$. However, methods applied in the proofs of theorems related to Besov and Triebel-Lizorkin spaces require a restriction of this range to $\alpha\in [-1/2,\infty)^d\setminus(-1/2,1/2)^d$. The same restraint has been used before (see for example \cite{NowakStempak}).

The Besov and Triebel-Li\-zor\-kin spaces were also studied in the setting of Hermite function expansions. For example in \cite{PetrushevXu} they were defined via norms
$$\Vert f\Vert_{F^{\alpha,q}_p}=\Big\Vert \Big(\sum\limits_{j=0}^{\infty}\left( 2^{\alpha j} \vert \phi_j* f(\cdot)\vert\right)^q\Big)^{1/q}\Big\Vert_p,$$
$$\Vert f\Vert_{B^{\alpha,q}_p}=\Big(\sum\limits_{j=0}^{\infty}\Big(\left\Vert  2^{\alpha j} \vert \phi_j* f(\cdot)\right\Vert_p\Big)^q\Big)^{1/q},$$
where the construction of functions $\phi_j$ was based on Hermite functions.

Many authors studied the above mentioned spaces introducing a lot of different definitions, often based on various decompositions (see for instance \cite{BuiPaluchTaibleson2, BuiPaluchTaibleson, Dziubanski, Epperson, Triebel}). Also the same spaces in Laguerre setting were frequently investigated (see \cite{BuiDuong2, BuiDuong3, KerPetPicXu}). In this paper we introduce definitions of Besov and Triebel-Li\-zor\-kin spaces in the setting of Laguerre expansions of Hermite type analogous to the definitions in \cite{BuiDuong}.

The main aim of the paper is to extend the result obtained in \cite{BuiDuong} for Hermite function expansions to Laguerre function expansions of Hermite type. We use the molecular decomposition analogous to \cite{BuiDuong}. We replace the Hermite operator $H=-\Delta+\vert x\vert^2$ by the operator 
$$\Lafi=-\Delta+\vert x\vert^2 +\sum\limits_{i=1}^d \frac{\alpha_i^2-1/4}{x_i^2},$$
whose eigenfunctions are the Laguerre functions of Hermite type. The results and methods in Sections 4 and 5 are similar to the corresponding ones in \cite{BuiDuong}.

The connection between the Hermite functions and the Laguerre functions of Hermite type is well known (see for example \cite{NowakStempak}). Therefore it seems natural to try to generalize the results obtained for the first functions to the latter.

We consider the Laguerre functions of Hermite type for the Laguerre type multi-index $\alpha\in (-1,\infty)^d$. In this paper the Laguerre functions of convolution type also appear but are used  only as a tool to characterize some Schwartz type spaces.

An important difference to \cite{BuiDuong} is the existence of the boundary of the considered domain which is $\mr^d_+$. This leads to problems nearby the hyperplanes orthogonal to the unit coordinate vectors. An important simplification comparing to \cite{BuiDuong} is the introduction of Schwartz spaces based on Theorem \ref{Schwartzcharacterization}. It allows us to avoid arduous calculations and replace them by more elegant reasoning.

The organization of this paper is as follows. In Section 2 we introduce the Laguerre functions of Hermite and convolution types and the operators associated with them. Then we define the Schwartz spaces. In the next section we study the operators associated with $\Lafi$. The final part of the section contains a fragment of the theory of subharmonic functions. In the last two section we demonstrate the results for Besov and Triebel-Li\-zor\-kin spaces.

We shall also make a frequent use, often without mentioning it in relevant places, of the two following facts: for any $A>0$ and $a\geq 0$, $\sup_{t>0}t^a\exp(-At)=C_{a,A}<\infty$; $\sum_{k\in\mathbb{N}^d\colon\vert k\vert= n}1={n+d-1\choose n}\lesssim (n+1)^{d-1}$ uniformly on $n\in\mathbb{N}$.

\subsection*{Notation}
Throughout the paper $\alpha=(\alpha_1,\ldots,\alpha_d)$ and $k=(k_1,\ldots,k_d)$ $\in\mathbb{N}^d$ are multi-indices, where $\mathbb{N}=\{0,1,\ldots\}$ and $d\geq1$ is the dimension. Unless stated otherwise $\alpha\in (-1,\infty)^d$, however in Section 4 and 5 we restrict it to $[-1/2,\infty)^d\setminus(-1/2,1/2)^d$. The length of multi-indices $k$ and $\alpha$ we denote by $\vert k\vert=k_1+\ldots+k_d$ and $\vert \alpha\vert=\alpha_1+\ldots+\alpha_d$. Note that $\vert \alpha\vert$ may be negative! Also we use the following notation: $\mathbb{N}_+=\{1,2,\ldots\}$ and $\mr^d_+=(0,\infty)^d$. The symbol $\lesssim$ we reserve for an inequality with some constant, independent of key parameters. We write $\simeq$ if there is $\lesssim$ and $\gtrsim$. The space $L^2(\mr^d_+,dx)$ will be simply denoted by $L^2(\mr^d_+)$. 
\subsection*{Acknowledgement}
The author would like to express his gratitude to Professor Krzysztof Stempak for indicating the topic and his constant support during the preparation of this paper.

\section{Laguerre Functions and Schwartz Spaces}
In this section we introduce the definitions of Laguerre functions of Hermite and convolution types and the operators associated with them. Our main aim is to define an analogue of tempered distributions related to appropriate test functions on $\mr^d_+$.
\subsection{Preliminaries}

We define the (countable) family of dyadic cubes $\mathcal{D}$ in $\mr^d_+$ as follows
$$\mathcal{D}=\bigcup\limits_{\nu\in\mathbb{Z}}\mathcal{D}_{\nu},$$
where
$$\mathcal{D}_{\nu}=\Big\{ \prod\limits_{j=1}^d \left( m_j 2^{\nu},(m_j+1)2^{\nu}\right]: m_1,\ldots,m_d\in\mathbb{N}\Big\}.$$
For a dyadic cube $Q$ we denote its length and center by $\ell(Q)$ and $x_Q$.
Immediately from the definition we have
$$\forall \nu\in\mathbb{Z}\qquad \mr_+^d=\bigcup_{Q\in\mathcal{D}_{\nu}}Q.$$

For $r>0$ we define the uncentered {\it maximal operator} $\mathcal{M}_r$ by
$$\mathcal{M}_rf(x)=\sup_{Q\ni x}\Big(\frac{1}{\vert Q\vert}\int\limits_Q \vert f(y)\vert^rdy\Big)^{1/r},\qquad x\in\mr^d.$$
$\mathcal{M}_r$ operates on measurable functions defined on $\mr^d$, where the supremum is taken over all cubes (not necessarily dyadic but with sides parallel to the axes of the coordinate system) $Q\subset\mr^d$ containing $x$.

The operator $\mathcal{M}_r$ is of weak type $(r,r)$ and of strong type $(p,p)$ for $p>r$. Moreover, we recall the Fefferman-Stein vector valued maximal inequality (see \cite{bigStein}): for $0<p<\infty$, $0<q\leq\infty$, $0<r<\min(p,q)$ and a sequence of functions $\{f_n\}_{n\in\mathbb{N}}$ we have
$$\Big\Vert\Big(\sum\limits_{n=0}^{\infty}\vert \mathcal{M}_rf_n\vert^q  \Big)^{1/q}\Big\Vert_p\lesssim \Big\Vert\Big(\sum\limits_{n=0}^{\infty}\vert f_n\vert^q  \Big)^{1/q}\Big\Vert_p.$$\\

The {\it Laguerre functions of Hermite type} of order $\alpha$ on $\mathbb{R}_+^d$ are the functions  
$$\varphi_k^{\alpha}(x)=\varphi_{k_1}^{\alpha_1}(x_1)\cdot\ldots\cdot\varphi_{k_d}^{\alpha_d}(x_d),\qquad x=(x_1,...,x_d)\in\mathbb{R}^d_+,$$
where $\varphi_{k_i}^{\alpha_i}$ is the one-dimensional Laguerre function of Hermite type of order $\alpha_i$
$$\varphi_{k_i}^{\alpha_i}(x_i)=\Big(\frac{2\Gamma(k_i+1)}{\Gamma(k_i+\alpha_i+1)}  \Big)^{1/2}L_{k_i}^{\alpha_i}(x_i^2)x_i^{\alpha_i+1/2}e^{-x_i^2/2},\qquad x_i>0,\ i=1,...,d,$$
whereas for $k\in\mathbb{N}$ and $\alpha>-1$, $L_k^{\alpha}$ denotes the Laguerre polynomial of degree $k$ and order $\alpha$ defined by
$$L_k^{\alpha}(x)=\frac{x^{-\alpha}e^x}{k!}\frac{d^k}{dx^k}\Big(e^{-x}x^{k+\alpha} \Big). $$
The functions $\{\varphi_k^{\alpha}:\ k\in\mathbb{N}^d\}$ form an orthonormal basis in $L^2(\mr_+^d,dx)$.

The {\it Laguerre functions of convolution type} of order $\alpha$ on $\mathbb{R}_+^d$ are the functions 
$$\ell_k^{\alpha}(x)=\ell_{k_1}^{\alpha_1}(x_1)\cdot\ldots\cdot\ell_{k_d}^{\alpha_d}(x_d),\qquad x=(x_1,...,x_d)\in\mathbb{R}^d_+,$$
where $\ell_{k_i}^{\alpha_i}$ is the one-dimensional Laguerre function of convolution type of order $\alpha_i$
$$\ell_{k_i}^{\alpha_i}(x_i)=\Big(\frac{2\Gamma(k_i+1)}{\Gamma(k_i+\alpha_i+1)}  \Big)^{1/2}L_{k_i}^{\alpha_i}(x_i^2)e^{-x_i^2/2},\qquad x_i>0,\ i=1,...,d.$$
Hence
\begin{equation}\label{phi&l}
\varphi_k^{\alpha}(x)=x^{\alpha+1/2}\ell_k^{\alpha}(x),
\end{equation}
where $x^{\alpha+1/2}=x_1^{\alpha_1+1/2}\cdot\ldots\cdot x_d^{\alpha_d+1/2}$.
The functions $\{\ell_k^{\alpha}:\ k\in\mathbb{N}^d\}$ form an orthonormal basis in $L^2(\mr_+^d,x^{2\alpha+1} dx)$.

We have the estimates (see for instance \cite[p.~99]{StempakTohoku} and \cite[p.~251]{Nowak})
\begin{equation}\label{unifbound}
\Vert \varphi_k^{\alpha}\Vert_{\infty}\lesssim 1,\quad \Vert \ell_k^{\alpha}\Vert_{\infty}\lesssim (1+\vert k\vert)^{c_{d,\alpha}},\qquad k\in\mathbb{N}^d,
\end{equation}
where $c_{d,\alpha}$ is a constant that depends only on the dimension and $\alpha$.

Consider the operators
\begin{equation}\label{definitionLafi}
\Lafi=-\Delta+V^{\alpha}(x),\qquad \Lal=-\Delta+\vert x\vert^2 +\sum\limits_{i=1}^d \frac{2\alpha_i+1}{x_i}\frac{\partial}{\partial x_i},
\end{equation}
where $V^{\alpha}(x)$ is the potential defined by
\begin{equation}\label{definitionofV}
V^{\alpha}(x)=\vert x\vert^2 +\sum\limits_{i=1}^d \frac{\alpha_i^2-1/4}{x_i^2}
\end{equation}

It is known that each $\varphi_k^{\alpha}$ is an eigenfuction of $\Lafi$ corresponding to the eigenvalue $4\vert k\vert +2\vert\alpha\vert +2d$. The eigenvalues depend only on the length of the multi-index so we use the notation 
$$\lambda_{\vert k\vert}^{\alpha}=4\vert k\vert+2\vert\alpha\vert+2d. $$
Similarly, each $\ell_k^{\alpha}$ is an eigenfuction of $\Lal$ corresponding to the eigenvalue $\lambda_{\vert k\vert}^{\alpha}$.
The symbol $\Lafi$ will also denote the self-adjoint extension of the operator defined above (initially considered on the domain $C_c^{\infty}$) in terms of eigenfunctions $\{ \varphi_k^{\alpha}\}$. We describe the construction in the next section.\\ 

\subsection{Schwartz spaces $\Safi$ and $\Sal$}
The following characterization of the Schwartz space $\mathcal{S}$ on $\mr^d$ in terms of the Hermite functions $\{h_k\}_{k\in\mathbb{N}^d}$ is known (see \cite[pp.~141-145]{ReedSimon}).

\begin{thm}\label{Schwartzcharacterization} Let $\phi\in L^2(\mr^d)$. Then $\phi\in\mathcal{S}(\mr^d)$ if and only if for every $N\in\mathbb{N}$
$$\langle \phi,h_k\rangle_{L^2(\mr^d,dx)}=O((1+\vert k\vert)^{-N}), $$
uniformly in $k\in\mathbb{N}^d$.\\
\end{thm}

Similarly, we define the Schwartz spaces $\Safi$ and $\Sal$ on $\mr_+^d$ associated with the Laguerre functions of Hermite and convolution type.
{\defin\label{defSafiSal}We define
$$\Safi=\left\{\phi\in L^2(\mr_+^d,dx)\colon\quad \forall N\in\mathbb{N^d}\quad \langle \phi,\varphi_k^{\alpha}\rangle_{L^2(\mr^d_+,dx)} =O\Big((1+\vert k\vert)^{-N}\Big),\quad k\in\mathbb{N}^d \right\},$$
$$\Sal=\Big\{\phi\in L^2(\mr_+^d,x^{2\alpha+1 }dx)\colon\ \forall N\in\mathbb{N^d}\ \langle \phi,\ell_k^{\alpha}\rangle_{L^2(\mr^d_+,x^{2\alpha +1}dx)} =O\Big((1+\vert k\vert)^{-N}\Big),\ k\in\mathbb{N}^d \Big\}.$$}
Henceforth we denote the inner product in $L^2(\mr_+^d,dx)$ simply by $\langle\cdot ,\cdot \rangle$.
 
Also, we need an analogue of the classical Schwartz space on $\mr^d_+$. A function $f$ defined on $\mr^d$ is called multi-even if
$$f(x_1,\ldots ,x_d)=f(\vert x_1\vert, \ldots,\vert x_d\vert),\qquad (x_1,\ldots,x_d)\in\mr^d.$$
We define Schwartz space $\Se$ on $\mr^d_+$, as the space of restrictions of multi-even functions from $\mathcal{S}(\mr^d)$ to $\mr_+^d$.

\begin{thm}\label{Saficharac}For $\alpha\in \left(-1,\infty\right)^d$ the following identities hold 
$$\Sal=\Se,\qquad\Safi=x^{\alpha+1/2}\cdot\Se.$$
\end{thm}
\noindent Note that when $\alpha=(-1/2,\ldots,-1/2)$, then $\Safi=\Se$.

We need the following Lemma.

{\lm\label{StempakLemma}If $f\in\Se(\mr_+^d)$, then $\frac{1}{x_i}\frac{\partial f}{\partial x_i}\in\Se(\mr^d_+)$, $i=1,\ldots,d.$}\\
For the proof see \cite[Prop.\ 2.3]{StempakExtension}.
{\cor\label{LafifinSe}If $f\in\Se(\mr^d_+)$, then $\Lal f\in\Se(\mr_+^d)$.}
\begin{proof}
Recall that
$$\Lal f(x)=\Big(-\Delta+\vert x\vert^2 -\sum\limits_{i=1}^d \frac{2\alpha_i+1}{x_i}\frac{\partial }{\partial x_i}\Big)f(x).$$
Obviously, $\big(-\Delta+\vert x\vert^2\big)$ maps $\Se$ into itself. Lemma \ref{StempakLemma} states that the remaining component also has this property.
\end{proof}

Now we present the proof of Theorem \ref{Saficharac}.

\begin{proof}
The second equality simply follows from the first equality and \eqref{phi&l}. Indeed, we have
$$\big\langle (\cdot)^{\alpha+1/2}\phi,\varphi_k^{\alpha}\big\rangle = \int\limits_{\mr^d_+}x^{\alpha+1/2}\phi(x)\varphi_k^{\alpha}(x)dx = \int\limits_{\mr^d_+}\phi(x)\ell_k^{\alpha}(x)x^{2\alpha+1}dx =   \left\langle \phi,\ell_k^{\alpha}\right\rangle_{L^2(\mr^d_+,x^{2\alpha +1}dx)}. $$
Hence, it suffices to prove the two inclusions associated to the first equality.

Firstly, we show the inclusion $\Sal\subset \Se$. We will prove it only in the one-dimensional case; it is easy to generalize it to the multidimensional case.\\
Let $f\in\Sal$, $\alpha>-1$. We show that $f$ can be extended to a Schwartz function on $\mr$. We have
\begin{equation}\label{equalinL2}
f=\sum\limits_{k=0}^{\infty}\langle f,\ell_k^{\alpha}\rangle_{L^2(\mr_+,x^{2\alpha +1}dx)} \ell_k^{\alpha}
\end{equation}
in $L^2(\mr_+,x^{2\alpha+1}dx)$, where $\langle\cdot,\cdot\rangle$ denotes the inner product in $L^2(\mr_+,x^{2\alpha+1}dx)$. We will show that the series on the right hand side of \eqref{equalinL2} is convergent for all $x\in\mr$ (we consider the functions $\ell_k^{\alpha}$ as naturally defined on $\mr$) and is a Schwartz function on $\mr$. Let
$$S_Mf=\sum\limits_{k=0}^{M}\langle f,\ell_k^{\alpha}\rangle_{L^2(\mr_+,x^{2\alpha +1}dx)} \ell_k^{\alpha}.$$
$S_Mf$ is a Cauchy sequence with respect to supremum norm on $\mr$. Indeed, using \eqref{unifbound} we obtain for any $N\in\mathbb{N}$
\begin{equation}\label{estim1}
\vert S_{M_1+M_2}f(x)-S_{M_1}f(x)\vert=\Big\vert\sum\limits_{k=M_1+1}^{M_1+M_2}\langle f,\ell_k^{\alpha}\rangle_{L^2(\mr_+,x^{2\alpha +1}dx)} \ell_k^{\alpha}(x)\Big\vert\lesssim \sum\limits_{k=M_1+1}^{M_1+M_2} (1+k)^{-N}, 
\end{equation}
hence we may make the expression above arbitrary small. The series is convergent in $L^2(\mr_+,x^{2\alpha+1}dx)$, so it has a subsequence convergent almost everywhere, thus the equality \eqref{equalinL2} holds also almost everywhere. Because of \eqref{estim1} and the continuity of $\ell_k^{\alpha}$ the function $f$ is continuous. Note that
$$(\ell_k^{\alpha}(x))'=-x\ell_k^{\alpha}(x)-2\frac{c_{k,\alpha}}{c_{k-1,\alpha+1}}x\ell^{\alpha+1}_{k-1}(x),$$
where $c_{k,\alpha}=\big(\frac{2\Gamma(k+1)}{\Gamma(k+\alpha+1)}  \big)^{1/2}$. Let $m\in\mathbb{N}$. Then there exists a constant $c$ such that
$$\big\vert \frac{d^m}{dx^m} \ell_k^{\alpha}(x)\big\vert\lesssim (1+k)^c $$
on $(\varepsilon,1/\varepsilon)$, for any $\varepsilon>0$. Hence, we may differentiate the series
$$\sum\limits_{k=0}^{\infty}\langle f,\ell_k^{\alpha}\rangle_{L^2(\mr_+,x^{2\alpha +1}dx)} \frac{d^m}{dx^m} \ell_k^{\alpha}(x)$$
term by term. Thus, $f$ is smooth and its every derivatives decrease rapidly. Also the function $f$ is even because the functions $\ell^{\alpha}_k$ have this property. 

Now, we prove the inclusion $\Se\subset \Sal$. We cannot restrict the reasoning to the one-dimensional case (see the proof of Lemma \ref{StempakLemma}).\\
Let $f\in\Se$. We will show that
$$\forall N\in\mathbb{N}\qquad \exists c>0\qquad \forall k\in\mathbb{N}^d\qquad \vert \langle f,\ell_k^{\alpha} \rangle_{L^2(\mr^d_+,x^{2\alpha +1}dx)}\vert\leq c(1+\vert k\vert)^{-N}.$$
By Corollary \ref{LafifinSe} it follows that $\Lal f\in\Se$. Recall that $\ell_k^{\alpha}=\frac{1}{\lambda_k^{\alpha}}\Lal \ell_k^{\alpha}$, where $\lambda_k^{\alpha}=4k+2\alpha+2$. Hence, for any $N\in\mathbb{N}_+$, by using the symmetry of $\Lafi$ $N$ times
$$
\vert\langle f,\ell_k^{\alpha}\rangle_{L^2(\mr^d_+,x^{2\alpha +1}dx)}\vert =\frac{1}{(\lambda_k^{\alpha})^N}\left\vert\langle (\Lal)^N f,\ell_k^{\alpha}\rangle_{L^2(\mr^d_+,x^{2\alpha +1}dx)}\right\vert\lesssim (1+\vert k\vert)^{-N}.$$
This finishes the proof. 
\end{proof}

\subsection{Further information of $\Safi$}
The space $\Safi$ is naturally determined by the family of `norms' $\{q^{\alpha}_N\}_{N\in \mathbb{N}}$ defined by 
$$q^{\alpha}_N(\phi)=\sup\limits_{k\in\mathbb{N}^d}(1+\vert k\vert)^N \vert\langle \phi,\varphi_k^{\alpha}\rangle\vert,   \qquad \phi\in L^2(\mr^d_+),$$
in the sense that
$$\Safi=\big\{ \phi\in L^2(\mr_+^d,dx)\colon\ \forall N\in\mathbb{N}\quad q_N^{\alpha}(\phi)<\infty\big\}.$$ 
We define also a family of `norms' $\{p^{\alpha}_r\}_{r>0}$ on $L^2(\mr_+^d)$ by
$$p_r^{\alpha}(\phi)=\sum\limits_{n=0}^{\infty}(n+1)^r \Big(\sum\limits_{\vert k\vert =n}\left\vert \langle\phi,\varphi_k^{\alpha}\rangle\right\vert^2\Big)^{1/2},\qquad \phi\in L^2(\mr_+^d). $$
Note that for $0<r_1<r_2$ and every $\phi\in L^2(\mr_+^d)$ there holds $p^{\alpha}_{r_1}(\phi)\leq p^{\alpha}_{r_2}(\phi)$.

{\prop\label{NormEquiv}The following equality is valid
$$\Safi=\big\{ \phi\in L^2(\mr_+^d,dx)\colon\ \forall r>0\quad p_r^{\alpha}(\phi)<\infty\big\}.$$}
\begin{proof}
It suffices to show that each $q_N^{\alpha}$ is dominated by a combination of  $p_r^{\alpha}$'s and vice versa. Let $N\in\mathbb{N}$ and $\phi\in\Safi$. We obtain
\begin{align*}
q_N^{\alpha}(\phi)=\sup\limits_{k\in\mathbb{N}^d}(1+\vert k\vert)^{N}\vert \langle\phi,\varphi_k^{\alpha}\rangle\vert &\leq \sum\limits_{k\in\mathbb{N}^d}(1+\vert k\vert)^{N}\vert \langle\phi,\varphi_k^{\alpha}\rangle\vert\\
&=\sum\limits_{n=0}^{\infty}(n+1)^N \sum\limits_{\vert k\vert =n}\left\vert \langle\phi,\varphi_k^{\alpha}\rangle\right\vert\\
&\lesssim\sum\limits_{n=0}^{\infty}(n+1)^{N+(d-1)/2} \Big(\sum\limits_{\vert k\vert =n}\left\vert \langle\phi,\varphi_k^{\alpha}\rangle\right\vert^2\Big)^{1/2}.
\end{align*}
Thus $q_N^{\alpha}(\phi)\lesssim p_{N+(d-1)/2}^{\alpha}(\phi)$. On the other hand, let $r>0$ and $\phi\in\Safi$. Then 
\begin{align*}
p_r^{\alpha}(\phi)=\sum\limits_{n=0}^{\infty}(1+n)^{(r+d+1)-(d+1)} \Big(\sum\limits_{\vert k\vert =n}\left\vert \langle\phi,\varphi_k^{\alpha}\rangle\right\vert^2\Big)^{1/2}&\leq q_{\lceil r\rceil +d+1}^{\alpha}(\phi)\sum\limits_{n=0}^{\infty} (n+1)^{-2}\lesssim  q_{\lceil r\rceil +d+1}^{\alpha}(\phi).
\end{align*}
\end{proof}

We introduce a topology in $\Safi$  defined by the family of norms $\{p^{\alpha}_r\}_{r>0}$ (see \cite[pp.~125-126]{ReedSimon}). The convergence $\phi_k\rightarrow \phi$ in $\Safi$ means that 
$$\forall r>0\qquad p_r^{\alpha}(\phi_k-\phi)\stackrel{k\rightarrow\infty}{\rightarrow} 0. $$
With this topology $\Safi$ becomes a Fr\'{e}chet space. The proof of Proposition \ref{NormEquiv} implies that the topologies defined on $\Safi$ by the families of norms $\{q^{\alpha}_N\}_{N\in \mathbb{N}}$ and $\{p^{\alpha}_r\}_{r>0}$ are the same (see 
\cite[pp.~131-132]{ReedSimon}).

We define $\Safid$ as the space of all linear continuous functionals on $\Safi$. Its elements we call {\it tempered distributions} on $\Safi$. The topology in $\Safid$ is the weakest topology for which every mapping 
$$\Safid\ni T\rightarrow \langle T,\phi\rangle,\qquad \phi\in\Safi, $$
is continuous. The convergence of a sequence of tempered distributions means the convergence for every test  function:
$$f_n\rightarrow f\quad\text{in }\Safid \Longleftrightarrow \forall\phi\in\Safi\quad \langle f_n,\phi\rangle \rightarrow \langle f,\phi\rangle. $$
The symbol $\langle\cdot,\cdot\rangle$ denotes the action of a distribution on a test function. The same symbol denotes the inner product in $L^2(\mr^d_+,dx)$ but we hope it won't lead to a confusion.

{\rem \label{ContinEquiv} A mapping $T:\Safi\rightarrow\mathbb{C}$ is continuous if and only if there exists $r>0$, that
$$\vert T(\phi)\vert \lesssim p_r^{\alpha}(\phi),\qquad \phi\in\Safi.$$
In particular, if $f\in\Safid$, then there exists $r>0$, such that  
$$\vert\langle f,\varphi_k^{\alpha}\rangle\vert \lesssim (1+\vert k\vert)^{r},\qquad k\in\mathbb{N}^d.$$
}

\section{Operators associated to  $\Lafi$}
In this section we study some families of operators based on the operator $\Lafi$. In particular, we consider the Poisson semigroup associated with $\Lafi$ and more general  families of operators which we use to define Besov and Triebel-Lizorkin spaces. Also, we formulate the Calder\'{o}n reproducing formula associated with $\Lafi$. At the end, we recall a part of the theory of subharmonic functions and present its application in this paper. 

\subsection{Definitions and basic facts}
The operator $\Lafi$ defined as in \eqref{definitionLafi} on the domain $C_c^{\infty}(\mr_+^d)$ is symmetric and non-negative. Its self-adjoint extension, which we denote by $\Lafi$ as well, is defined on the domain  
$$\mathfrak{D}(\Lafi)=\Big\{\phi\in L^2(\mr^d_+):\quad \sum\limits_{n=0}^{\infty}\sum\limits_{\vert k\vert=n}\left\vert\lambda^{\alpha}_n\langle \phi,\varphi_k^{\alpha}\rangle\right\vert^2<\infty\Big\}   $$
by
$$\Lafi \phi=\sum\limits_{n=0}^{\infty}\lambda^{\alpha}_n\sum\limits_{\vert k\vert=n}\langle \phi,\varphi_k^{\alpha}\rangle\varphi_k^{\alpha},$$
and the convergence is in $L^2(\mr^d_+)$. The self-adjoint extension of $\Lafi$ is non-negative and its spectrum is the discrete set
$$\sigma(\Lafi)=\{\lambda_n^{\alpha}\colon n\in\mathbb{N}\}.$$

For a function $\mathfrak{m}$ defined on $\sigma(\Lafi)$, according to the spectral theorem, we consider the operator $\mathfrak{m}(\Lafi)$
\begin{equation}\label{spectraldef}
\mathfrak{m}(\Lafi)\phi=\sum\limits_{n=0}^{\infty}\mathfrak{m}(\lambda^{\alpha}_n)\sum\limits_{\vert k\vert=n}\langle \phi,\varphi^{\alpha}_k\rangle\varphi^{\alpha}_k,\qquad \phi\in L^2(\mr_+^d),
\end{equation}
which is defined on the domain
\begin{equation}\label{domainofm(L)}
\mathfrak{D}\big( \mathfrak{m}(\Lafi)\big)=\Big\{\phi\in L^2(\mr^d_+)\colon \sum\limits_{n=0}^{\infty}\sum\limits_{\vert k\vert=n}\left\vert\mathfrak{m}(\lambda^{\alpha}_n)\langle \phi,\varphi_k^{\alpha}\rangle\right\vert^2<\infty\Big\}.
\end{equation}
If $\mathfrak{m}$ is bounded on $\sigma(\Lafi)$, then $\mathfrak{D}\big( \mathfrak{m}(\Lafi)\big)=L^2(\mr^d_+)$ and $\mathfrak{m}(\Lafi)$ is bounded on $L^2(\mr^d_+)$, and if $\mathfrak{m}$ is real-valued, then $\mathfrak{m}(\Lafi)$ is self-adjoint.\\

We consider the {\it Poisson semigroup} associated to $\Lafi$, $\big\{e^{-t\sqrt{\Lafi}}\big\}_{t> 0}$, and more generally, for $m\in\mathbb{N}$, the family of operators
$$\Big\{\Big(t\sqrt{\Lafi}\Big)^me^{-t\sqrt{\Lafi}}\Big\}_{t>0}.$$
We use the notation
$$P^{\alpha}_{t,m}=\Big(t\sqrt{\Lafi}\Big)^me^{-t\sqrt{\Lafi}}.$$
According to \eqref{spectraldef} the operators $\Patm$ are spectrally defined on $L^2(\mr^d_+)$ by
\begin{equation}\label{spectraldefPatm}
\Patm \phi=\sum\limits_{n=0}^{\infty}\big(t\sqrt{\lambda^{\alpha}_n}\big)^m e^{-t\sqrt{\lambda^{\alpha}_n}}\sum\limits_{\vert k\vert=n} \langle \phi,\varphi_k^{\alpha}\rangle\varphi_k^{\alpha}.
\end{equation}
When $m=0$ we write $\Pat$.  

Moreover, the operators $\Patm$ are integral operators:
\begin{equation}\label{defPatm}
\forall \phi\in L^2(\mr_+^d)\qquad  \Patm \phi(x)= \int\limits_{\mr^d_+}\patm(x,y)\phi(y)dy,\qquad x\in\mr_+^d,\ t>0,
\end{equation}
and their kernels are defined by
\begin{equation}\label{kerneldefPatm}
\patm (x,y)=\sum\limits_{n=0}^{\infty}\big(t\sqrt{\lambda^{\alpha}_n}\big)^m e^{-t\sqrt{\lambda^{\alpha}_n}}\sum\limits_{k=\vert n\vert}\varphi_k^{\alpha}(x)\varphi_k^{\alpha}(y),\qquad x,y\in\mr^d_+.
\end{equation}
When $m=0$ we write $\pat$ for the Poisson kernel.

{\prop\label{patminSafi} Let $t>0$ and $m\in\mathbb{N}$.
\begin{enumerate}[(i)]
\item If $\phi_j\in \Safi$ converges to $\phi\in\Safi$ in $\Safi$, then $\phi_j$ also converges to $\phi$ in $ L^2(\mr_+^d)$. 
\item The series in \eqref{spectraldefPatm} converges to $\Patm \phi$ in $\Safi$ (for any $\phi\in L^2(\mr^d_+)$).
\item For any fixed $y\in\mr^d_+$ the series in \eqref{kerneldefPatm} converges to $\patm(\cdot,y)$ in $\Safi$.
\end{enumerate}}
\begin{proof} We begin with (i). For the simplicity assume that $\phi_j\rightarrow0$ in $\Safi$. Then
\begin{align*}
\Vert \phi_j\Vert_2=\Big(\sum\limits_{k\in\mathbb{N}^d}\vert\langle\phi_j,\varphi_k\rangle\vert^2 \Big)^{1/2}\leq \sum\limits_{n=0}^{\infty}\big( \sum\limits_{\vert k\vert=n}\vert \langle\phi_j,\varphi_k\rangle\vert^2\big)^{1/2}
\end{align*}
and the last quantity converges to $0$ because $\phi_j$ converges to $0$ in $\Safi$.

Now we pass to the proof of (ii). Let $r>0$. For $t,m$ and $\phi$ as in the claim and positive integer $M$ we compute

\begin{align*}
p_r^{\alpha}\Big(\sum\limits_{n=M}^{\infty}\big(t\sqrt{\lambda^{\alpha}_n}\big)^m e^{-t\sqrt{\lambda^{\alpha}_n}}\sum\limits_{\vert k\vert=n} \langle \phi,\varphi_k^{\alpha}\rangle\varphi_k^{\alpha}   \Big)&=\sum\limits_{n=M}^{\infty}(1+n)^r\Big(\sum\limits_{\vert k\vert=n}\big\vert \big( t\sqrt{\lambda_{\vert k\vert}^{\alpha}}\big)^{m}e^{-t\sqrt{\lambda_{\vert k\vert}^{\alpha}}}\langle\phi,\varphi_k^{\alpha}\rangle \big\vert^2\Big)^{\frac{1}{2}}\\
&=  \sum\limits_{n=M}^{\infty}(1+n)^r \big( t\sqrt{\lambda_{n}^{\alpha}}\big)^{m}e^{-t\sqrt{\lambda_{n}^{\alpha}}}   \Big(\sum\limits_{\vert k\vert=n}\left\vert \langle\phi,\varphi_k^{\alpha}\rangle \right\vert^2\Big)^{\frac{1}{2}}\\
&\lesssim \sum\limits_{n=M}^{\infty}(1+n)^{r+(d-1)/2} \big( t\sqrt{\lambda_{n}^{\alpha}}\big)^{m}e^{-t\sqrt{\lambda_{n}^{\alpha}}}\Vert\phi\Vert_2.
\end{align*}
The series is convergent so it tends to $0$ when $M\rightarrow \infty$.\\
The proof of (iii) follows the similar argument along with the uniform boundedness of the functions $\varphi_k^{\alpha}$. 
\end{proof}

\subsection{Kernel estimates}
In this paper we will not consider the heat semigroup associated with $\Lafi$, that is the semigroup of operators $\left\{ e^{-t\Lafi}\right\}_{t> 0}$, in the context of Besov and Triebel-Lizorkin spaces. However, it will be useful in kernel estimates.

The heat semigroup is a semigroup of integral operators on $L^2(\mr_+^d)$:
$$e^{-t\Lafi}\phi(x)=\int\limits_{\mr^d_+}G_t^{\alpha}(x,y)\phi(y)dy,\qquad \phi\in L^2(\mr_+^d),\ x\in\mr^d_+,\ t>0, $$
where
$$G_t^{\alpha}(x,y)=\sum\limits_{n=0}^{\infty} e^{-t\lambda^{\alpha}_n}\sum\limits_{k=\vert n\vert}\varphi_k^{\alpha}(x)\varphi_k^{\alpha}(y),\qquad x,y\in\mr^d_+,$$
and there is an explicit formula
$$G_t^{\alpha}(x,y)=(\sinh(2t))^{-d}\exp\Big(-\frac{1}{2}\coth(2t)(\vert x\vert^2+\vert y\vert^2)\Big)\prod\limits_{i=1}^d\sqrt{x_iy_i}\mathcal{I}_{\alpha_i}\Big(\frac{x_iy_i}{\sinh(2t)}\Big).$$
Here $\mathcal{I}_{\nu}$, $\nu\in\mr$, denotes the modified Bessel function of first kind and order $\nu$, which is smooth and positive on $(0,\infty)$.

Before we pass to the estimates of $\patm$ we give some technical lemmas. 

{\lm\label{heatkernelestim}For $\alpha\in [-1/2,\infty)^d$ the following estimate holds
$$ G^{\alpha}_t(x,y)\lesssim \left\{\begin{array}{ll} 
e^{-dt}e^{-\frac{\vert x-y\vert^2}{2}},         &t>1,\\
t^{-d/2} e^{-\frac{\vert x-y\vert^2}{4t}}, &0<t\leq 1,\\
\end{array}\right.\qquad x,y\in\mr_+^d.$$
}For the proof see \cite[Lemma 2.1]{StempakWrobel} and the comment following it.

The subordination identity
$$e^{-t\sqrt{a}}=\frac{1}{4\sqrt{\pi}}\int\limits_0^{\infty}te^{-\frac{t^2}{4u}}e^{-ua}\frac{du}{u^{3/2}},\quad a>0, $$
and a computation show that $\pat(x,y)$ and $G_t^{\alpha}(x,y)$ are related by
\begin{equation}\label{patandG}
\pat(x,y)=\frac{1}{4\sqrt{\pi}}\int\limits_0^{\infty}te^{-\frac{t^2}{4u}}G_u^{\alpha}(x,y)\frac{du}{u^{3/2}}.
\end{equation}
We shall need the following representation of $\patm(x,y)$.

{\lm\label{patmestimates} Let $t>0$, $x,y\in\mr^d_+$, $m\in\mathbb{N}$ and $\alpha\in [-1/2,\infty)^d$. Then
$$\patm(x,y)=\frac{(-1)^{m+1}t^m}{2\sqrt{\pi}}\int\limits_0^{\infty}\partial_t^{m+1}\big( e^{-\frac{t^2}{4u}}\big)G_u^{\alpha}(x,y)\frac{du}{u^{1/2}}.$$}
\begin{proof}

Firstly we shall check that
\begin{equation}\label{difpatm}
\patm(x,y)=(-t)^m\partial^m_t\pat(x,y).
\end{equation}
It suffices to use \eqref{kerneldefPatm} and justify the possibility of differentiating the series on the right hand side of \eqref{kerneldefPatm} taken with $m=0$, term by term successively $m$ times with respect to $t$, for any fixed $x,y\in\mr^d_+$. This is indeed possible due to an appropriate application of Lebesgue's dominated convergence theorem.

Now, using \eqref{patandG} and \eqref{difpatm} and entering with $\partial_t^m$ under the integral sign (this is again legitimate due to Lebesgue's dominated convergence theorem) we finally arrive to
$$\patm(x,y)=\frac{(-t)^m}{4\sqrt{\pi}}\int\limits_0^{\infty}\partial_t^{m}\big( te^{-\frac{t^2}{4u}}\big)G_u^{\alpha}(x,y)\frac{du}{u^{3/2}}=\frac{(-1)^{m+1}t^m}{2\sqrt{\pi}}\int\limits_0^{\infty}\partial_t^{m+1}\big(e^{-\frac{t^2}{4u}}\big)G_u^{\alpha}(x,y)\frac{du}{u^{1/2}}.$$
\end{proof}

{\prop\label{estimpatm}For every $m\in\mathbb{N}_+$ and $\alpha\in [-1/2,\infty)^d$ we have 
$$\vert\patm(x,y)\vert\lesssim\frac{t^m}{(t+\vert x-y\vert)^{d+m}},\qquad t>0,\quad x,y\in\mr_+^d.$$}
\begin{proof}
Note that for every $m\in\mathbb{N}_+$ we have
$$\Big\vert \partial_t^m e^{-\frac{t^2}{4u}}\Big\vert \lesssim e^{-\frac{t^2}{8u}}u^{-m/2},\qquad u>0,\quad t>0.$$
For the proof see \cite[p.~410]{BuiDuong}. Thus Lemma \ref{patmestimates} implies
$$\vert \patm (x,y)\vert \lesssim t^m \int\limits_0^{\infty} \Big\vert \partial_t^{m+1}\big(e^{-\frac{t^2}{4u}}\big)\Big\vert   G_u^{\alpha}(x,y)\frac{du}{u^{1/2}}\lesssim t^m \int\limits_0^{\infty} e^{-\frac{t^2}{8u}}u^{-(m+1)/2}   G_u^{\alpha}(x,y) \frac{du}{u^{1/2}}.$$
We now split the integration over $(0,\infty)$ onto the intervals $(0,1)$ and $(1,\infty)$, denote the resulting integrals by $I_0$ and $I_{\infty}$, and estimate them. For $I_0$ we write using the substitution $u=(t^2+\vert x-y\vert^2)/8w$. Lemma \ref{heatkernelestim} yields
\begin{align*}
I_0\lesssim \int\limits_0^{1}u^{-(m+d+2)/2}e^{-\frac{t^2+\vert x-y\vert^2}{8u}}du &\lesssim \int\limits_{\frac{t^2+\vert x-y\vert^2}{8}}^{\infty}e^{-w} \frac{w^{(m+d-2)/2}}{(t^2+\vert x-y\vert^2)^{(m+d)/2}}dw\\
&\lesssim \Gamma\Big(\frac{m+d}{2}\Big) \frac{1}{(t+\vert x-y\vert)^{m+d}},
\end{align*}
and for $I_{\infty}$ we estimate
$$I_{\infty}\lesssim \int\limits_1^{\infty} e^{-ud}u^{-(m+2)/2} e^{-\frac{t^2+\vert x-y\vert^2}{8u}}du\lesssim \frac{1}{\left( t+\vert x-y\vert\right)^{m+d}} \int\limits_1^{\infty} e^{-ud} u^{(d-2)/2}du\lesssim \frac{1}{\left( t+\vert x-y\vert\right)^{m+d}}.$$
The required estimate in Proposition \ref{estimpatm} follows.
\end{proof}

\subsection{$\Patm$ on tempered distributions}
Thanks to Proposition \ref{patminSafi}, for $f\in\Safid$, $t>0$ and $m\in\mathbb{N}$, we may define $\Patm f$ as a function on $\mr^d_+$ by
\begin{equation}\label{defdistrib}
\Patm f(x)=\langle f, \patm(x,\cdot)\rangle,\qquad x\in\mr^d_+. 
\end{equation}

{\prop\label{PatmfinSafi} Let $f,t$ and $m$ be as above. Then $\Patm f\in\Safi$. Moreover, the mapping $T_f:\Safi\rightarrow\mathbb{C}$ defined as
$$T_f(\phi)=\int\limits_{\mr^d_+} \Patm f(x)\phi(x)dx,\qquad \phi\in\Safi,$$
is a tempered distribution.}
\begin{proof}
Let $f\in\Safid$ and $r>0$ be such that $\vert\langle f,\varphi_k^{\alpha}\rangle\vert\lesssim (1+\vert k\vert)^r$, $k\in\mathbb{N}^d$ (see Remark \ref{ContinEquiv}). Firstly, we check that $\Patm f\in L^2(\mr_+^d)$. Using Proposition \ref{patminSafi} and Minkowski's inequality we compute
\begin{align*}
\Vert \Patm f\Vert_2=\Big( \int\limits_{\mr_+^d} \big\vert \langle f,\patm (x,\cdot)\rangle\big\vert^2 dx\Big)^{1/2}&=\Big( \int\limits_{\mr_+^d} \Big\vert \Big\langle f,\sum\limits_{n=0}^{\infty} \big(t\sqrt{\lambda_n^{\alpha}}\big)^m e^{-t\sqrt{\lambda_n^{\alpha}}}\sum\limits_{\vert k\vert=n}\varphi_k^{\alpha}(x)\varphi_k^{\alpha}   \Big\rangle\Big\vert^2 dx\Big)^{1/2}\\
&= \Big( \int\limits_{\mr_+^d} \Big\vert \sum\limits_{n=0}^{\infty} \big(t\sqrt{\lambda_n^{\alpha}}\big)^m e^{-t\sqrt{\lambda_n^{\alpha}}}\sum\limits_{\vert k\vert=n}\varphi_k^{\alpha}(x) \langle f,\varphi_k^{\alpha}\rangle \Big\vert^2 dx\Big)^{1/2}\\
&\lesssim \sum\limits_{n=0}^{\infty} \big(t\sqrt{\lambda_n^{\alpha}}\big)^m e^{-t\sqrt{\lambda_n^{\alpha}}} (1+n)^r  \sum\limits_{\vert k\vert=n} \Vert \varphi_k^{\alpha}\Vert_2\\
&\lesssim \sum\limits_{n=0}^{\infty} \big(t\sqrt{\lambda_n^{\alpha}}\big)^m e^{-t\sqrt{\lambda_n^{\alpha}}} (1+n)^{r+d-1}
\end{align*}
and the last sum is finite.

Now we shall prove that $\Patm f\in\Safi$. Fix $j\in\mathbb{N}^d$. Proposition \ref{patminSafi} implies that for any $N\in\mathbb{N}$, using Fubini's theorem, we have
\begin{align*}
\langle \Patm f,\varphi_j^{\alpha}\rangle &=\int\limits_{\mr^d_+}\Big\langle f,\sum\limits_{n=0}^{\infty} \big(t\sqrt{\lambda_n^{\alpha}}\big)^m e^{-t\sqrt{\lambda_n^{\alpha}}}\sum\limits_{\vert k\vert=n}\varphi_k^{\alpha}(x)\varphi_k^{\alpha}\Big\rangle \varphi_j^{\alpha}(x)dx\\
&=\int\limits_{\mr^d_+}\sum\limits_{n=0}^{\infty} \big(t\sqrt{\lambda_n^{\alpha}}\big)^m e^{-t\sqrt{\lambda_n^{\alpha}}}\sum\limits_{\vert k\vert=n}   \langle f,\varphi_k^{\alpha}\rangle\varphi_k^{\alpha}(x)\varphi_j^{\alpha}(x)dx\\
&=\big(t\sqrt{\lambda_{\vert j\vert}^{\alpha}}\big)^m e^{-t\sqrt{\lambda_{\vert j\vert}^{\alpha}}}\langle f,\varphi_j^{\alpha}\rangle\\
&=O\big((1+\vert j\vert)^{-N}\big).
\end{align*}
Now, we will show that $T_f\in\Safid$. It suffices to prove that there exists $\tilde{r}>0$ that $\vert T_f(\phi)\vert \lesssim p_{\tilde{r}}^{\alpha}(\phi),\ \phi\in\Safi$. Let $\phi\in\Safi$. Using Proposition \ref{patminSafi}, Fubini's theorem and the Cauchy-Schwarz inequality we obtain 
\begin{align*}
\vert T_f(\phi)\vert&=\Big\vert \int\limits_{\mr^d_+}\Big\langle f,\sum\limits_{n=0}^{\infty} \big(t\sqrt{\lambda_n^{\alpha}}\big)^m e^{-t\sqrt{\lambda_n^{\alpha}}}\sum\limits_{\vert k\vert=n}\varphi_k^{\alpha}(x)\varphi_k^{\alpha}\Big\rangle \phi(x)dx\Big\vert\\
&\leq  \sum\limits_{n=0}^{\infty} \big(t\sqrt{\lambda_n^{\alpha}}\big)^m e^{-t\sqrt{\lambda_n^{\alpha}}}\sum\limits_{\vert k\vert=n} \Big\vert \int\limits_{\mr^d_+}\left\langle f,\varphi_k^{\alpha}\right\rangle \varphi_k^{\alpha}(x) \phi(x) dx\Big\vert\\
&\lesssim \sum\limits_{n=0}^{\infty} \big(t\sqrt{\lambda_n^{\alpha}}\big)^m e^{-t\sqrt{\lambda_n^{\alpha}}}\sum\limits_{\vert k\vert=n} (1+n)^{r} \left\vert\langle \phi,\varphi_k^{\alpha}\rangle\right\vert\\
&\lesssim \sum\limits_{n=0}^{\infty} \big(t\sqrt{\lambda_n^{\alpha}}\big)^m e^{-t\sqrt{\lambda_n^{\alpha}}}(1+n)^{r+d-1}\Big(\sum\limits_{\vert k\vert=n}  \big\vert\langle \phi,\varphi_k^{\alpha}\rangle\big\vert^2\Big)^{1/2}\\
&\lesssim p_{r+d-1}^{\alpha}(\phi).
\end{align*}
\end{proof}

{\prop The operator $\Patm$, $t>0,\ m\in\mathbb{N}$, maps $\Safi$ into itself continuously.}
\begin{proof}
Let $\phi\in\Safi$. By Proposition \ref{patminSafi} (i), 

$\Patm\phi\in\Safi$. For $\phi_j\rightarrow0$ in $\Safi$ and fixed $r>0$ we compute
\begin{align*}
p_r^{\alpha}( \Patm\phi_j)&=\sum\limits_{n=0}^{\infty}(n+1)^r\Big( \sum\limits_{\vert k\vert=n} \big\vert \big\langle  \big(t\sqrt{\Lafi}\big)^m e^{-t\sqrt{\Lafi}}\phi_j,\varphi_k^{\alpha}\big\rangle\big\vert^2\Big)^{1/2}\\
&=\sum\limits_{n=0}^{\infty}(n+1)^r\Big( \sum\limits_{\vert k\vert=n} \big\vert \big(t\sqrt{\lambda_{\vert k\vert}^{\alpha}}\big)^{m}e^{-t\sqrt{\lambda_{\vert k\vert}^{\alpha}}} \left\langle  \phi_j,\varphi_{k}^{\alpha}\right\rangle\big\vert^2\Big)^{1/2}\\
&=\sum\limits_{n=0}^{\infty}(n+1)^r \big(t\sqrt{\lambda_{n}^{\alpha}}\big)^{m}e^{-t\sqrt{\lambda_{n}^{\alpha}}}\Big( \sum\limits_{\vert k\vert=n} \vert  \langle  \phi_j,\varphi_{k}^{\alpha}\rangle\vert^2\Big)^{1/2}\\
&\lesssim p_r^{\alpha}(\phi_j).
\end{align*}
Because of $p_r^{\alpha}(\phi_j)\rightarrow0$, $j\rightarrow\infty$, we have $p_r^{\alpha}\left( \Patm\phi_j\right)\rightarrow0$, $j\rightarrow\infty$.
\end{proof}
{\rem\label{remtempdistr} If $f\in\Safid$, then we also may define $\Patm f$ as tempered distribution by
$$\langle\Patm f,\phi\rangle:=\langle f,\Patm \phi\rangle,\qquad \phi\in\Safi,$$
and this definition coincides with the definition of $\Patm f$ as a function.}
Indeed, for $\phi\in\Safi$ applying Proposition \ref{patminSafi} and Fubini's theorem (a verification of the appropriate assumption will be given later) we obtain
\begin{align*}
\int\limits_{\mr^d_+}\Patm f(x)\phi(x)dx&=\int\limits_{\mr^d_+}\Big\langle f, \sum\limits_{n=0}^{\infty} \big(t\sqrt{\lambda_n^{\alpha}}\big)^m e^{-t\sqrt{\lambda_n^{\alpha}}}\sum\limits_{\vert k\vert=n}\varphi_k^{\alpha}(x)\varphi_k^{\alpha} \Big\rangle\phi(x)  dx\\
&=\sum\limits_{n=0}^{\infty} \big(t\sqrt{\lambda_n^{\alpha}}\big)^m e^{-t\sqrt{\lambda_n^{\alpha}}}\sum\limits_{\vert k\vert=n} \langle f,\varphi_k^{\alpha}\rangle\int\limits_{\mr^d_+}\varphi_k^{\alpha}(x)\phi(x)  dx\\
&= \Big\langle f,\sum\limits_{n=0}^{\infty} \big(t\sqrt{\lambda_n^{\alpha}}\big)^m e^{-t\sqrt{\lambda_n^{\alpha}}}\sum\limits_{\vert k\vert=n}\varphi_k^{\alpha} \int\limits_{\mr^d_+}\phi(x)\varphi_k^{\alpha}(x) dx\Big\rangle\\
&= \langle f,\Patm\phi  \rangle.
\end{align*}
Now, we verify the assumption of Fubini's theorem. For some $r>0$, the Cauchy-Schwarz inequality yields
\begin{align*}
\sum\limits_{n=0}^{\infty} \big(t\sqrt{\lambda_n^{\alpha}}\big)^m e^{-t\sqrt{\lambda_n^{\alpha}}}\sum\limits_{\vert k\vert=n} \vert\langle f,\varphi_k^{\alpha}\rangle\vert\int\limits_{\mr^d_+}\vert\varphi_k^{\alpha}(x)\phi(x)\vert  dx&\leq \sum\limits_{n=0}^{\infty} \big(t\sqrt{\lambda_n^{\alpha}}\big)^m e^{-t\sqrt{\lambda_n^{\alpha}}}\sum\limits_{\vert k\vert=n} (1+n)^r\Vert\phi\Vert_2\\
&\lesssim \sum\limits_{n=0}^{\infty} \big(t\sqrt{\lambda_n^{\alpha}}\big)^m e^{-t\sqrt{\lambda_n^{\alpha}}} (1+n)^{r+d-1}\Vert\phi\Vert_2\\
&\lesssim \Vert\phi\Vert_2.
\end{align*}

\subsection{Calder\'{o}n reproducing formula}
{\prop Let $f\in\Safid$. Then
\begin{itemize}
\item[(i)] for all $m\in\mathbb{N}$, $\lim\limits_{t\rightarrow\infty}\Patm f=0$ in $\Safid$;
\item[(ii)] for all $m\in\mathbb{N}_+$, $\lim\limits_{t\rightarrow 0^+}\Patm f=0$ in $\Safid$;
\item[(iii)] $\lim\limits_{t\rightarrow 0^+}(I-P^{\alpha}_t)f=0$ in $\Safid$.
\end{itemize}}\label{limits}
\begin{proof}
(i) Definition \ref{remtempdistr} implies that it suffices to prove that for any $\phi\in\Safi$ and $r>0$ there is
$$\lim\limits_{t\rightarrow \infty}p^{\alpha}_r( \Patm \phi)=0. $$
Hence, applying Fubini's theorem, we compute
\begin{align*}
p_r^{\alpha}( \Patm\phi)&=\sum\limits_{n=0}^{\infty}(n+1)^r\Big(\sum\limits_{\vert k\vert=n}\vert \langle\Patm\phi,\varphi_k^{\alpha}\rangle\vert^2\Big)^{1/2}\\
&=\sum\limits_{n=0}^{\infty}(n+1)^r\Big(\sum\limits_{\vert k\vert=n}\Big\vert \sum\limits_{i=0}^{\infty}\big(t\sqrt{\lambda_i^{\alpha}}\big)^m e^{-t\sqrt{\lambda_i^{\alpha}}}	\sum\limits_{\vert j\vert=i}\langle\phi,\varphi_j^{\alpha}\rangle\langle\varphi_j^{\alpha},\varphi_k^{\alpha}\rangle\Big\vert^2\Big)^{1/2}\\
&=\sum\limits_{n=0}^{\infty}(n+1)^r\Big(\sum\limits_{\vert k\vert=n}\Big\vert\big(t\sqrt{\lambda_{\vert k\vert}^{\alpha}}\big)^me^{-t\sqrt{\lambda_{\vert k\vert}^{\alpha}}}\langle\phi,\varphi_k^{\alpha}\rangle\Big\vert^2\Big)^{1/2}\\
&=\sum\limits_{n=0}^{\infty}(n+1)^r \big(t\sqrt{\lambda_{n}^{\alpha}}\big)^m e^{-t\sqrt{\lambda_{n}^{\alpha}}}\Big(\sum\limits_{\vert k\vert=n}\vert\langle\phi,\varphi_k^{\alpha}\rangle\vert^2\Big)^{1/2},
\end{align*}
and the last quantity converges to $0$ when $t\rightarrow\infty$, because  $\phi\in\Safi$.\\
$\ $\\
(ii) As in (i), one just have to replace $t\rightarrow\infty$ by $t\rightarrow 0^+$ in the last step.\\
$\ $\\
(iii) Similarly to (i), we will show that for any $\phi\in\Safi$ and $r>0$ we have
$$\lim\limits_{t\rightarrow 0^+}p_r^{\alpha}( (I-P^{\alpha}_t)\phi)=0.$$
Thus
\begin{align*}
p_r^{\alpha}( (I-P^{\alpha}_t)\phi) &=\sum\limits_{n=0}^{\infty}(n+1)^r\Big(\sum\limits_{\vert k\vert=n}\Big\vert \Big\langle\phi-\sum\limits_{i=0}^{\infty}e^{-t\sqrt{\lambda_i^{\alpha}}}\sum\limits_{\vert j\vert=i}\langle\phi,\varphi_j^{\alpha}\rangle\varphi_j^{\alpha},\varphi_k^{\alpha}\Big\rangle\Big\vert^2\Big)^{1/2}\\
&=\sum\limits_{n=0}^{\infty}(n+1)^r\Big(\sum\limits_{\vert k\vert=n}\Big\vert \langle\phi,\varphi_k^{\alpha}\rangle-e^{-t\sqrt{\lambda_{\vert k\vert}^{\alpha}}}\langle\phi,\varphi_k^{\alpha}\rangle\Big\vert^2\Big)^{1/2}\\
&= \sum\limits_{n=0}^{\infty}(n+1)^r\Big(1-e^{-t\sqrt{\lambda_n^{\alpha}}}\Big)\Big(\sum\limits_{\vert k\vert=n}\vert \langle\phi,\varphi_k^{\alpha}\rangle\vert^2\Big)^{1/2},
\end{align*}
and the last expression converges to $0$ when $t\rightarrow 0^+$ because $\phi\in\Safi$.
\end{proof}

{\rem In the following proposition we deal with integration of a one-parameter family of distributions. We understand it as follows
$$\Big\langle \int_0^{\infty} \Patm f \frac{dt}{t},\phi\Big\rangle=\Big\langle f,\int_0^{\infty} \Patm\phi \frac{dt}{t}\Big\rangle,\qquad \phi\in\Safi. $$}

{\prop[Calder\'{o}n reproducing formula] Let $m_1,m_2\in\mathbb{N}_+$, $m=m_1+m_2$ and $f\in\Safid$. Then
$$f=\frac{2^m}{(m-1)!}\int\limits_0^{\infty}P^{\alpha}_{t,m_1}P^{\alpha}_{t,m_2}f\frac{dt}{t} $$
in $\Safid$.}{\label{Calformula}}
\begin{proof}
Since the operators $P^{\alpha}_{t,m_i}$, $i=1,2$, commute the equality above may be written as
$$f=\frac{1}{(m-1)!}\int\limits_0^{\infty}P^{\alpha}_{2t,m} f\frac{dt}{t}.$$
Let $0<\delta<\tau<\infty$. We integrate by parts
\begin{align*}
2^{-m}\int\limits_{\delta}^{\tau}P^{\alpha}_{2t,m}f\frac{dt}{t}&=-\left(t\sqrt{\Lafi}\right)^{m-1}\frac{1}{2}e^{-2t\sqrt{\Lafi}}f\Big\vert_{\delta}^{\tau}+\frac{m-1}{2}\int\limits_{\delta}^{\tau}\left(t\sqrt{\Lafi}\right)^{m-2}\sqrt{\Lafi} e^{-2t\sqrt{\Lafi}}fdt\\
&=-\left(t\sqrt{\Lafi}\right)^{m-1}\frac{1}{2}e^{-2t\sqrt{\Lafi}}f\Big\vert_{\delta}^{\tau}-\left(t\sqrt{\Lafi}\right)^{m-2}\frac{m-1}{4}e^{-2t\sqrt{\Lafi}}f\Big\vert_{\delta}^{\tau}\\
&\qquad\qquad\qquad+\frac{(m-1)(m-2)}{4}\int\limits_{\delta}^{\tau}\left(t\sqrt{\Lafi}\right)^{m-3}\sqrt{\Lafi} e^{-2t\sqrt{\Lafi}}fdt\\
&=\ldots=-\sum\limits_{n=0}^{m-1}\frac{(m-1)!}{2^{m-n}n!}\left(t\sqrt{\Lafi}\right)^ne^{-2t\sqrt{\Lafi}}f\Big\vert_{\delta}^{\tau}\\
&=-\sum\limits_{n=0}^{m-1}\frac{(m-1)!}{2^{m-n}n!}\left[\left(\tau\sqrt{\Lafi}\right)^{n}e^{-\sqrt{2\tau\Lafi}}-\left(\delta\sqrt{\Lafi}\right)^{n}e^{-\sqrt{2\delta\Lafi}}     \right]f.
\end{align*}
Proposition \ref{limits} implies
$$\lim\limits_{\delta\rightarrow 0^+}\left(\delta\sqrt{\Lafi}\right)^n e^{-2\delta\sqrt{\Lafi}}f=0$$
in $\Safid$ for $n\in\mathbb{N}_+$, 
$$\lim\limits_{\delta\rightarrow 0^+} e^{-2\delta\sqrt{\Lafi}}f=f$$
in $\Safid$ and
$$\lim\limits_{\tau\rightarrow \infty}\left(\tau\sqrt{\Lafi}\right)^n e^{-2\tau\sqrt{\Lafi}}f=0$$
in $\Safid$ for $n\in\mathbb{N}$.\\
Letting $\delta\rightarrow 0^+$ and $\tau\rightarrow\infty$ in the result above we obtain
$$\int\limits_0^{\infty}P^{\alpha}_{t,m_1}P^{\alpha}_{t,m_2}f\frac{dt}{t}=\frac{(m-1)!}{2^{m}}f$$
in $\Safid$ and this finishes the proof.
\end{proof}
\subsection{Subharmonic functions}
Let $\Omega\subset\mr^d$ be a nonempty open set. A locally integrable in $\Omega$ real-valued function $g$ is called subharmonic if $\Delta g\geq 0$ in the sense of $\mathcal{D}'(\Omega)$, where $\mathcal{D}'(\Omega)$ denotes the space of (classical) distributions on $\Omega$ and $\mathcal{D}(\Omega)=C^{\infty}_c(\Omega)$. If $g\in C^2(\Omega)$, then the condition is equivalent to the non-negativity of $\Delta g$ as a function. 

Throughout this subsection, for $f\in\Safid$ and $m\in\mathbb{N}$, we will consider the function
$$\Pam f(x)=t^{-m}\Patm f(x) $$
as a function of the variables $(x,t)\in \mr^d_+\times\mr_+$. Note that \eqref{kerneldefPatm}, \eqref{defdistrib} and Proposition \ref{patminSafi} yield 
\begin{equation}\label{pamformula}
\Pam f(x)=\sum\limits_{n=0}^{\infty}\big(\sqrt{\lambda^{\alpha}_n}\big)^m e^{-t\sqrt{\lambda^{\alpha}_n}}\sum\limits_{k=\vert n\vert}\varphi_k^{\alpha}(x)\langle f,\varphi_k^{\alpha}\rangle.
\end{equation}

{\lm\label{PamCinfty} $\Pam f(x)\in C^{\infty}(\mr^d_+\times\mr_+)$ and we can differentiate (with respect to $x$ or $t$) the series \eqref{pamformula} term by term.}
\begin{proof}
The equality
$$\frac{\partial}{\partial x_i} \varphi_k^{\alpha}(x)=-2k_i^{1/2}\varphi_{k-e_i}^{\alpha+e_i}+\big(\frac{2\alpha_i-1}{2x_i}-x_i\big)\varphi_k^{\alpha}(x) $$
and \eqref{unifbound} imply that $\vert\frac{\partial}{\partial x_i}\varphi_k^{\alpha}\vert\lesssim (k_i+1)^{1/2}$ on $(\varepsilon,1/\varepsilon)^d$, where $e_i$ is the $i$-th coordinate vector in $\mr^d$, $\varepsilon$ any positive number and, by convention, $\varphi_{k-e_i}^{\alpha+e_i}=0$ if $k_i=0$. Moreover, $\vert\frac{d}{d t}e^{-t\sqrt{\lambda_n^{\alpha}}}\vert\lesssim \sqrt{\lambda_n^{\alpha}}$ on $(\varepsilon,\infty)$. Also we apply Remark \ref{ContinEquiv}. Therefore we can differentiate (with respect to $x$ or $t$) the series in \eqref{pamformula} term by term. Similarly for higher derivatives.
\end{proof}

{\lm\label{Pamsubharmonic} $\vert \Pam f(x)\vert^2$ is a subharmonic function in $\mr_+^d\times\mr_+$.}
\begin{proof}
Firstly, we will check that
\begin{equation}\label{laplacianPam}
\big(\Delta+\frac{\partial^2}{\partial t^2}\big)\Pam f(x) =V^{\alpha}(x)\Pam f(x),
\end{equation}
where $\Delta$ denotes the Laplacian on $\mr_+^d$ and $V^{\alpha}(x)$ is defined in \eqref{definitionofV}. Applying Lemma \ref{PamCinfty} we obtain
\begin{align*}
\big(\Delta+\frac{\partial^2}{\partial t^2}\big)\Pam f(x)&=\sum\limits_{n=0}^{\infty}\big(\sqrt{\lambda^{\alpha}_n}\big)^m \sum\limits_{k=\vert n\vert}\langle f,\varphi_k^{\alpha}\rangle \big(\Delta+\frac{\partial^2}{\partial t^2}\big)e^{-t\sqrt{\lambda^{\alpha}_n}}\varphi_k^{\alpha}(x)\\
&=\sum\limits_{n=0}^{\infty}\big(\sqrt{\lambda^{\alpha}_n}\big)^m \sum\limits_{k=\vert n\vert}\langle f,\varphi_k^{\alpha}\rangle e^{-t\sqrt{\lambda^{\alpha}_n}}\big(\Delta\varphi_k^{\alpha}(x)+\lambda^{\alpha}_n\varphi_k^{\alpha}(x)\big).
\end{align*}
Using the fact that $\varphi^{\alpha}_k$ are the eigenfunctions of $\Lafi$ completes the proof of \eqref{laplacianPam}. Now we shall prove the claim. We denote the gradient in $\mr^d_+\times\mr_+$ by $\nabla$ and compute
$$\big(\Delta+\frac{\partial^2}{\partial t^2}\big)\vert \Pam f(x)\vert^2= 2\Vert\nabla\Pam f(x)\Vert^2 +2 V^{\alpha}(x)\vert \Pam f(x)\vert^2.$$
The obtained quantity is non-negative, thus the proof is completed.
\end{proof}

{\prop\label{Patmsubharm} For $\alpha\in [-1/2,\infty)^d\setminus(-1/2,1/2)^d$ let $g(x,t)$ denote the multi-even extension (with respect to $x$) of $\Pam f(x)$ in $\mr^d\times\mr_+$. Then $\vert g(x,t)\vert^2$ is a subharmonic function in $\mr^d\times\mr_+$.}
\begin{proof}
Note that Lemma \ref{Pamsubharmonic} says that $\vert g(x,t)\vert^2$ is subharmonic in $\mr^d_+\times\mr_+$. We shall prove that it is also subharmonic in $\mr^d\times\mr_+$.

Remark that Theorem \ref{Saficharac}, Proposition \ref{PatmfinSafi} and Lemma \ref{PamCinfty} imply that $g(x,t)\in C^1(\mr^d\times\mr_+)$ (not necessarily in $C^2$; in case $\alpha_i=-1/2$ for some $i=1,\ldots,d$, $g(x,t)$ for fixed the rest of variables belongs to $\mathcal{S}(\mr)$), $g_{x_ix_i}(x,t)$, $i=1,\ldots,d,$ is defined almost everywhere and is locally integrable in $\mr^d\times\mr_+$, whereas $g_{tt}(x,t)$ is continuous in $\mr^d\times\mr_+$.

Hence, we shall check the appropriate condition in $\mathcal{D}'(\mr^d\times\mr_+)$. For the sake of clarity we present the computation in case $d=1$; the case of higher dimensions follows easily. Let $\phi(x,t)\in C^{\infty}_c(\mr\times\mr_+)$. Firstly, note that
\begin{align*}
&\int\limits_{-\infty}^{\infty} g(x,t)\phi_{xx}(x,t)dx\\
&= \big(g(0^-,t)-g(0^+,t)\big)\phi_x(0,t)-\big(g_x(0^-,t)-g_x(0^+,t)\big)\phi(0,t)+\int\limits_{-\infty}^{\infty} g_{xx}(x,t)\phi(x,t)dx\\
&=\int\limits_{-\infty}^{\infty} g_{xx}(x,t)\phi(x,t)dx.
\end{align*}
Also,
$$\int\limits_{0}^{\infty} g(x,t)\phi_{tt}(x,t)dt=\int\limits_{0}^{\infty} g_{tt}(x,t)\phi(x,t)dt.$$
Thus, applying Fubini's theorem we obtain
\begin{align*}
\int\limits_0^{\infty}\int\limits_{-\infty}^{\infty} g(x,t)\big(\phi_{xx}(x,t)+\phi_{tt}(x,t)\big)dxdt=\int\limits_0^{\infty}\int\limits_{-\infty}^{\infty} \big(g_{xx}(x,t)+g_{tt}(x,t)\big)\phi(x,t)dxdt.
\end{align*} 
Moreover, the functions $g_{xx}(x,t)$ and $g_{tt}(x,t)$ are even (with respect to $x$). Hence, Lemma \ref{Pamsubharmonic} yields that $\Delta_{(x,t)}g\geq 0$ in $\mathcal{D}'(\mr\times\mr_+)$. The proof is completed.
\end{proof}

{\lm\label{subharmlemma}Let $g$ be a non-negative subharmonic function in $\Omega\subset\mr^d$. Then for every cube $Q$ such that $\overline{2Q}\subset\Omega$, $1<\mu\leq 2$, $r>0$, and almost every $x\in Q$ the following inequality holds
$$ g(x)\leq \Big(\frac{1}{\vert \mu Q\vert}\int\limits_{\mu Q}  g(y)^r dy\Big)^{1/r},$$
where $\mu Q$ is a cube with the same center as $Q$ and the length $\ell(\mu Q)=\mu\ell(Q)$.}\\
For the proof see \cite[pp.~10-11]{AuscherBenAli}.

\section{Homogeneous Besov spaces associated to $\Lafi$}
In this section we introduce the homogeneous Besov spaces in the setting of the Laguerre function expansions of Hermite type and present the main theorems. The proofs are given in the latter part of this section. From now on we assume $\alpha\in [-1/2,\infty)^d\setminus(-1/2,1/2)^d$.
\subsection{Definitions and theorems}
In what follows given $0<p,q\leq\infty$ we denote $r_0=\min(1,p,q)$ and then, given $d\in\mathbb{N}_+$ and $\sigma\in\mr$, we define $m_0=d+\max(\sigma,0)+\lfloor d(1/r_0-1)\rfloor +1$.
{\defin For $\sigma\in\mr$, $0<p,q\leq\infty$ and $m\in\mathbb{N}$, such that $m>m_0$ we define the {\it homogeneous Besov spaces} $\Bifi$ by
$$\Bifi=\Big\{ f\in\Safid:\quad \Vert f\Vert_{\Bifi}=\Big(\int\limits_0^{\infty} \left( t^{-\sigma}\left\Vert \Patm f\right\Vert_p\right)^q\frac{dt}{t}\Big)^{1/q}<\infty\Big\}.$$}

Our main tool will be the molecular decomposition analogous to the one introduced in \cite{BuiDuong}. In the following definition we use the fact that $\mathfrak{D}\big((\sqrt{\Lafi})^{j_0}\big)\subset\mathfrak{D}\big((\sqrt{\Lafi})^{j}\big)$ for every $0<j\leq j_0$ (see \eqref{domainofm(L)}).

{\defin Let $0<p\leq\infty$, $\sigma\in\mr$, and $M,N\in\mathbb{N}_+$. A function $a\in L^2(\mr^d_+)$ is called $(\Lafi, M,N,\sigma,p)$ {\it molecule} if there exist a function $b\in\mathfrak{D}\big(( \sqrt{\Lafi})^{2M}\big)$, $\nu\in\mathbb{Z}$, and a dyadic cube $Q\in\mathcal{D}_{\nu}$, such that
\begin{enumerate}[(i)]
\item $\left( \sqrt{\Lafi}\right)^M b=a$;
\item $\big\vert \left(\sqrt{\Lafi}\right)^j b(x)\big\vert \lesssim 2^{\nu(M-j+\sigma)}\vert Q\vert^{-1/p}\big( 1+\frac{\vert x-x_Q\vert}{2^{\nu}}\big)^{-d-N}$  for $j=0,\ldots, 2M$ and $x\in\mr^d_+$ a.e.
\end{enumerate}
}

{\thm\label{moldecBesov}Let $\sigma\in\mr$ and $0<p,q\leq\infty$.
\begin{enumerate}[(a)]
\item For $M,N\in\mathbb{N}_+$ and $m\in\mathbb{N}$, such that $m>m_0$, if $f\in\Bifi$, then there exist a sequence of $(\Lafi, M,N,\sigma,p)$ molecules $\{a_Q\}_{Q\in\mathcal{D}}$ and a sequence of coefficients $\{s_Q\}_{Q\in\mathcal{D}}$ such that
$$f=\sum\limits_{Q\in\mathcal{D}}s_Qa_Q $$
in $\Safid$. Moreover,
\begin{equation}\label{simcoeffBesov}
\Big[ \sum\limits_{\nu\in\mathbb{Z}}\Big(\sum\limits_{Q\in\mathcal{D}_{\nu}}\vert s_Q\vert^p\Big)^{q/p}\Big]^{1/q}\simeq \Vert f\Vert_{\Bifi}
\end{equation}
\item Conversely, if $M>\max(d/r_0-\sigma,m)$, $N>d(1/r_0-1)$ and $m>\max(\sigma,0)+N+d$, and for a sequence of $(\Lafi, M,N,\sigma,p)$ molecules $\{a_Q\}_{Q\in\mathcal{D}}$ and a sequence of coefficients  $\{s_Q\}_{Q\in\mathcal{D}}$ satisfying
$$\Big[ \sum\limits_{\nu\in\mathbb{Z}}\Big(\sum\limits_{Q\in\mathcal{D}_{\nu}}\vert s_Q\vert^p\Big)^{q/p}\Big]^{1/q}<\infty,$$
the series $\sum\limits_{Q\in\mathcal{D}}s_Qa_Q $ converges in $\Safid$, then its sum $f$ is in $\Bifi$. Moreover,
$$\Vert f\Vert_{\Bifi}\lesssim \Big[ \sum\limits_{\nu\in\mathbb{Z}}\Big(\sum\limits_{Q\in\mathcal{D}_{\nu}}\vert s_Q\vert^p\Big)^{q/p}\Big]^{1/q}.$$
\end{enumerate}}

{\thm\label{BesovEquiv}Let $\sigma\in\mr$, $0<p,q\leq\infty$, $m_1,m_2\in\mathbb{N}$ and $m_1,m_2>m_0$. Then the spaces $\dot{B}^{\sigma,\Lafi ,m_1}_{p,q}$ and $\dot{B}^{\sigma,\Lafi ,m_2}_{p,q}$ coincide and their norms are equivalent.}\\

Theorem \ref{BesovEquiv} implies that the Besov space is independent of the $m$ index provided $m>m_0$ hence we drop the index in further part of the paper and denote the Besov space by $\dot{B}^{\sigma,\Lafi}_{p,q}$.\\

We also state the embedding theorem for Besov spaces $\dot{B}^{\sigma,\Lafi}_{p,q}$.
{\thm\label{BesovEmbed}Let $\sigma\in\mr$, $0<p,q\leq\infty$. Then
\begin{enumerate}[(i)]
\item if $0<q_1\leq q_2 $, then $\dot{B}^{\sigma,\Lafi}_{p,q_1}\hookrightarrow \dot{B}^{\sigma,\Lafi}_{p,q_2}$,
\item if $\sigma_1\geq\sigma_2$ and $\sigma_1-d/p_1=\sigma_2-d/p_2$, then $\dot{B}^{\sigma_1,\Lafi}_{p_1,q}\hookrightarrow \dot{B}^{\sigma_2,\Lafi}_{p_2,q}$.
\end{enumerate}}$\ $

To prove the above theorems we shall need the following two lemmas.
{\lm\label{molest}Suppose that $\sigma\in\mr$ and the assumptions on $M,N,m$ and $\{a_Q\}_{Q\in\mathcal{D}}$ are as in Theorem \ref{moldecBesov} (b). Then for $\nu\in\mathbb{Z}$ and $Q\in\mathcal{D}_{\nu}$ we have
\begin{enumerate}[(i)]
\item for $t\leq 2^{\nu}$
$$\vert \Patm a_Q(x)  \vert\lesssim \vert Q\vert^{-1/p} 2^{\nu\sigma}\Big( \frac{t}{2^{\nu}}\Big)^{m-N-d}\Big(1+\frac{\vert x-x_Q\vert}{2^{\nu}}\Big)^{-N-d},$$
uniformly with respect to $Q$, $x$, $t$, $\sigma$ and $\nu$;
\item for $t>2^{\nu}$
$$\vert \Patm a_Q(x)  \vert\lesssim \vert Q\vert^{-1/p} 2^{\nu\sigma}\Big( \frac{2^{\nu}}{t}\Big)^{M}\Big(1+\frac{\vert x-x_Q\vert}{t}\Big)^{-N-d}$$
uniformly with respect to $Q$, $x$, $t$, $\sigma$ and $\nu$.
\end{enumerate}
}$\ $\\
\begin{proof}
(i) Proposition \ref{estimpatm}, the definition of $(\Lafi, M,N,\sigma,p)$ molecules and the formula $a_Q=\big(\sqrt{\Lafi}\big)^M b_Q$, together with \eqref{defPatm} and Proposition \ref{estimpatm} imply for $x\in\mr^d_+$
\begin{align*}
\vert \Patm a_Q(x) \vert& =\Big\vert t^{m-N} P^{\alpha}_{t,N}\Big(\Big(\sqrt{\Lafi}\Big)^{m-N}a_Q  \Big)(x)\Big\vert\\
&\lesssim \int\limits_{\mr^d_+} \frac{t^N}{(t+\vert x-y\vert)^{N+d}}t^{m-N}\Big\vert  \Big(\Big(\sqrt{\Lafi}\Big)^{m-N+M}b_Q  \Big)(y)  \Big\vert dy\\
&\lesssim \int\limits_{\mr^d_+} \frac{t^{m-N-d}}{\Big(1+\frac{\vert x-y\vert}{t}\Big)^{(N+d)}} 2^{\nu(N-m+\sigma)}\vert Q\vert^{-1/p}\Big(1+\frac{\vert y-x_Q\vert}{2^{\nu}}  \Big)^{-(d+N)} dy   \\
&\lesssim 2^{\nu\sigma}\Big(\frac{t}{2^{\nu}} \Big)^{m-N}\vert Q\vert^{-1/p}\int\limits_{\mr^d} \frac{t^{-d}}{\Big(1+\frac{\vert x-y\vert}{t}\Big)^{(N+d)}}    \Big(1+\frac{\vert y-x_Q\vert}{2^{\nu}}  \Big)^{-(d+N)} dy.
\end{align*}
We now split the integration over $\mr^d$ onto the sets $\{y\in\mr^d\colon  \vert x-y\vert \leq \vert x-x_Q\vert/2 \}$ and $\{y\in\mr^d\colon  \vert x-y\vert > \vert x-x_Q\vert/2 \}$, denote the resulting integrals by $I_1(x)$ and $I_2(x)$, and estimate them.

For $I_1$ note that if $\vert x-y\vert\leq\frac{1}{2}\vert x-x_Q\vert$, then $\vert y-x_Q\vert \sim\vert x-x_Q\vert$.
 
Hence
$$I_1(x) \lesssim \Big(1+\frac{\vert x-x_Q\vert}{2^{\nu}}  \Big)^{-(d+N)}\int\limits_{\mr^d}t^{-d} \Big(1+\frac{\vert x-y\vert}{t}\Big)^{-(d+N)} dy\lesssim    \Big(1+\frac{\vert x-x_Q\vert}{2^{\nu}}  \Big)^{-(d+N)}.$$
Estimating $I_2$ note that for $y$ such that $\vert x-y\vert >\frac{1}{2}\vert x-x_Q\vert$ it holds
$$\Big(1+\frac{\vert x-y\vert}{t}\Big)^{-(d+N)}\lesssim  \Big(1+\frac{\vert x-x_Q\vert}{t}\Big)^{-(d+N)}.$$
Thus
\begin{align*}
I_2(x)&\lesssim \int\limits_{\{y\in\mr^d\colon  \vert x-y\vert > \vert x-x_Q\vert/2 \}} \Big(1+\frac{\vert x-y\vert}{t}\Big)^{-(d+N)} t^{-d}  \Big(1+\frac{\vert y-x_Q\vert}{2^{\nu}}  \Big)^{-(d+N)} dy\\
&\lesssim  \Big(\frac{t}{2^{\nu}} \Big)^{-d} \Big(1+\frac{\vert x-x_Q\vert}{t}\Big)^{-(d+N)}\int\limits_{\mr^d} 2^{-d\nu}\Big(1+\frac{\vert y-x_Q\vert}{2^{\nu}}  \Big)^{-(d+N)} dy\\
&\lesssim \Big(\frac{t}{2^{\nu}} \Big)^{-d}\Big(1+\frac{\vert x-x_Q\vert}{2^{\nu}}\Big)^{-(d+N)},
\end{align*}
and this completes the proof of (i).\\
(ii) Putting  $a_Q=\left(\sqrt{\Lafi}\right)^M b_Q$, for $x\in\mr^d_+$ we have
\begin{align*}
\vert \Patm a_Q(x)\vert&=\vert t^{-M} P^{\alpha}_{t,m+M} b_Q(x)\vert\\
&\lesssim t^{-M}\int\limits_{\mr_+^d} t^{-d} \Big(1+\frac{\vert x-y\vert}{t}\Big)^{-(m+M+d)} \vert b_Q(y)\vert dy\\
&\lesssim t^{-M} \int\limits_{\mr_+^d} t^{-d} \Big(1+\frac{\vert x-y\vert}{t}\Big)^{-(N+d)}2^{\nu(M+\sigma)}\vert Q\vert^{-1/p}\Big(1+\frac{\vert y-x_Q\vert}{2^{\nu}}  \Big)^{-d-N} dy \\
&\lesssim \vert Q\vert^{-1/p}  2^{\nu\sigma}\Big(\frac{2^{\nu}}{t} \Big)^M\int\limits_{\mr_+^d} t^{-d} \Big[\Big(1+\frac{\vert x-y\vert}{t}\Big) \Big(1+\frac{\vert y-x_Q\vert}{2^{\nu}}  \Big)\Big]^{-d-N} dy. 
\end{align*} 
Arguing similarly as in (i) we get
$$\vert \Patm a_Q(x)  \vert\lesssim \vert Q\vert^{-1/p} 2^{\nu\sigma}\Big( \frac{2^{\nu}}{t}\Big)^{M}\Big(1+\frac{\vert x-x_Q\vert}{t}\Big)^{-N-d},$$
and this finishes the proof.
\end{proof}

{\lm\label{FrazJawLemma}Let $N>0,\ \nu,\mu\in\mathbb{Z}$, $\nu\leq\mu$, and $\left\{f_Q\right\}_{Q\in\mathcal{D}_{\nu}}$ be a sequence of functions satisfying
$$\vert f_Q(x)\vert\lesssim\Big( 1+\frac{\vert x-x_Q\vert}{2^{\mu}}\Big)^{-(d+N)}$$
uniformly with respect to $Q$, $x\in\mr^d_+$, $\nu$ and $\mu$. Then for $r>\frac{d}{d+N}$ and a sequence of coefficients $\left\{\tilde{s}_Q\right\}_{Q\in\mathcal{D}_{\nu}}$ it is true that
$$\sum\limits_{Q\in\mathcal{D}_{\nu}} \vert \tilde{s}_Q\vert \vert f_Q(x)\vert\lesssim 2^{(\mu-\nu)d/r}\mathcal{M}_r\Big(\sum\limits_{Q\in\mathcal{D}_{\nu}} \vert \tilde{s}_Q\vert \chi_Q\Big)(x) $$
uniformly with respect to $x\in\mr^d_+$, $\nu$ and $\mu$.}\\
For the proof see \cite[pp.~147-148]{FrazierJawerth}.

\subsection{Proofs of Theorems \ref{moldecBesov}, \ref{BesovEquiv}, \ref{BesovEmbed}}
\begin{proof}[Proof of Theorem \ref{moldecBesov}] We begin with the proof of (a). The reproducing Calder\'{o}n formula (Proposition \ref{Calformula}) implies that for $f\in\Bifi$
$$f=c_{m,M,N}\int\limits_0^{\infty} P^{\alpha}_{t,M+N}\Patm f\frac{dt}{t}$$
in $\Safi$, where $c_{m,M,N}=\frac{2^{m+M+N}}{(m+M+N-1)!}$. Hence
$$f=c_{m,M,N}\sum\limits_{\nu\in\mathbb{Z}}\int\limits_{2^{\nu}}^{2^{\nu+1}}P^{\alpha}_{t,M+N}\Patm f\frac{dt}{t}=c_{m,M,N}\sum\limits_{\nu\in\mathbb{Z}}\sum\limits_{Q\in\mathcal{D}_{\nu}}\int\limits_{2^{\nu}}^{2^{\nu+1}}P^{\alpha}_{t,M+N}(\Patm f\cdot\chi_Q)\frac{dt}{t}$$
in $\Safid$. The dyadic cube decomposition is possible because all of the integrals above we may consider as functions (see Definition \ref{remtempdistr}). Now, for $\nu\in\mathbb{Z}$ and $Q\in\mathcal{D}_{\nu}$ we define
$$s_Q=2^{-\nu\sigma}\vert Q\vert^{1/p}\sup\limits_{(y,t)\in Q\times (2^{\nu},2^{\nu+1}]} \vert \Patm f(y)\vert$$
and $a_Q=\left(\sqrt{\Lafi}\right)^M b_Q$, where
$$b_Q(x)=\frac{c_{m,M,N}}{s_Q}\int\limits_{2^{\nu}}^{2^{\nu+1}}t^M P^{\alpha}_{t,N} ( \Patm f\cdot\chi_Q)(x)\frac{dt}{t}.$$
Obviously, $f=\sum\limits_{Q\in\mathcal{D}}s_Qa_Q$ in $\Safid$. Moreover, using Proposition \ref{estimpatm} one can check that $a_Q$ are $(\Lafi,M,N,\sigma,p)$ molecules. 

It remains to prove \eqref{simcoeffBesov}. Firstly, we compute
\begin{align*}
\Vert f\Vert_{\Bifi}&=\Big[\sum\limits_{\nu\in\mathbb{Z}}\int\limits_{2^{\nu}}^{2^{\nu+1}}\Big(t^{-\sigma p}\sum\limits_{Q\in\mathcal{D_{\nu}}}\int\limits_Q \vert \Patm f(y)\vert^p dy\Big)^{q/p}\frac{dt}{t} \Big]^{1/q}\\
&\lesssim\Big[\sum\limits_{\nu\in\mathbb{Z}}\int\limits_{2^{\nu}}^{2^{\nu+1}}\Big(\sum\limits_{Q\in\mathcal{D_{\nu}}}  \vert s_Q\vert^p 2^{\nu p\sigma}t^{-\sigma p} \Big)^{q/p}\frac{dt}{t} \Big]^{1/q}\\
&\lesssim \Big[ \sum\limits_{\nu\in\mathbb{Z}}\Big(\sum\limits_{Q\in\mathcal{D}_{\nu}}\vert s_Q\vert^p\Big)^{q/p}\Big]^{1/q}.
\end{align*}
To prove the opposite inequality in \eqref{simcoeffBesov} note that
\begin{equation}\label{FormulaFors_Q}
\Big\Vert \sum\limits_{Q\in\mathcal{D}_{\nu}}\vert Q\vert^{-1/p}\vert s_Q\vert \chi_Q \Big\Vert_p=\Big(\sum\limits_{Q\in\mathcal{D}_{\nu}}\vert s_Q\vert^p\Big)^{1/p}.
\end{equation}
Fix $0<r<r_0$ and $Q\in\mathcal{D}_{\nu}$. Then $ \tilde{Q}=Q\times (2^{\nu},2^{\nu+1}]$ is a cube in $\mr^d\times\mr_+$ and $\overline{2Q}\subset\mr^d\times\mr_+$. Treating the function $\Patm f(x)$ as a function of variables $(x,t)$ defined in $\mr^d\times\mr_+$ as in Proposition \ref{Patmsubharm} and applying the theorem and Lemma \ref{subharmlemma} we obtain
$$\sup_{(y,t)\in  \tilde{Q}}\Big\vert \Big(\sqrt{\Lafi}\Big)^m e^{-t\sqrt{\Lafi}}f(y)\Big\vert\lesssim \Big(\frac{1}{\vert \frac{3}{2}\tilde{Q}\vert}\iint\limits_{ \frac{3}{2}\tilde{Q}} \Big\vert \Big(\sqrt{\Lafi}\Big)^m e^{-t\sqrt{\Lafi}}f(x)\Big\vert^r dxdt\Big)^{1/r}.$$
Note that there is $\vert\tilde{Q}\vert\sim 2^{\nu}\vert Q\vert$ and $t\sim 2^{\nu}$ for $t\in \tilde{Q}$. Hence for any $x\in Q$ 
\begin{align*}
\sup_{(y,t)\in  \tilde{Q}}\vert \Patm f(y)\vert&\lesssim \Big(\int\limits_{\frac{3}{4}2^{\nu}}^{\frac{9}{8}2^{\nu+1}}\frac{1}{\vert \frac{3}{2}Q\vert}\int\limits_{ \frac{3}{2}Q} \vert \Patm f(y)\vert^r dy\frac{dt}{t}\Big)^{1/r}\lesssim \Big(\int\limits_{\frac{3}{4}2^{\nu}}^{\frac{9}{8}2^{\nu+1}} [\mathcal{M}_r (\Patm f)(x)]^r\frac{dt}{t}\Big)^{1/r}.
\end{align*}
Thus 
$$\sum\limits_{Q\in\mathcal{D}_{\nu}}\vert Q\vert^{-1/p}\vert s_Q\vert \chi_Q(x)\lesssim 2^{-\nu\sigma}\Big(\int\limits_{\frac{3}{4}2^{\nu}}^{\frac{9}{8}2^{\nu+1}} [\mathcal{M}_r (\Patm f)(x)]^r\frac{dt}{t}\Big)^{1/r}.$$
The results above give
$$\Big(\sum\limits_{Q\in\mathcal{D}_{\nu}}\vert s_Q\vert^p\Big)^{1/p}\lesssim \Big\Vert 2^{-\nu\sigma} \Big(\int\limits_{\frac{3}{4}2^{\nu}}^{\frac{9}{8}2^{\nu+1}} [\mathcal{M}_r (\Patm f)]^r\frac{dt}{t}\Big)^{1/r}   \Big\Vert_p\lesssim 2^{-\nu\sigma}\Big\Vert \int\limits_{\frac{3}{4}2^{\nu}}^{\frac{9}{8}2^{\nu+1}} [\mathcal{M}_r (\Patm f)]^r\frac{dt}{t}   \Big\Vert_{p/r}^{1/r}.$$
Since $p/r>1$, the Minkowski inequality yields
$$\Big\Vert \int\limits_{\frac{3}{4}2^{\nu}}^{\frac{9}{8}2^{\nu+1}} [\mathcal{M}_r (\Patm f)]^r\frac{dt}{t}   \Big\Vert_{p/r}\leq \int\limits_{\frac{3}{4}2^{\nu}}^{\frac{9}{8}2^{\nu+1}} \Vert  [\mathcal{M}_r (\Patm f)]^r   \Vert_{p/r} \frac{dt}{t},$$
hence, applying strong type $(p,p)$ of the operator $\mathcal{M}_r$ and H\"{o}lder's inequality we obtain
\begin{align*}
\Big(\sum\limits_{Q\in\mathcal{D}_{\nu}}\vert s_Q\vert^p\Big)^{1/p}&\lesssim 2^{-\nu\sigma} \Big[ \int\limits_{\frac{3}{4}2^{\nu}}^{\frac{9}{8}2^{\nu+1}} \Vert  \Patm f \Vert_{p}^r \frac{dt}{t}\Big]^{1/r}\lesssim  \Big[ \int\limits_{\frac{3}{4}2^{\nu}}^{\frac{9}{8}2^{\nu+1}} \big( t^{-\sigma}\Vert  \Patm f \Vert_{p}\big)^q \frac{dt}{t}\Big]^{1/q}.
\end{align*}
Finally,
$$\Big[\sum\limits_{\nu\in\mathbb{Z}}\Big(\sum\limits_{Q\in\mathcal{D}_{\nu}}\vert s_Q\vert^p\Big)^{q/p}\Big]^{1/q}\lesssim \Big[\sum\limits_{\nu\in\mathbb{Z}}\int\limits_{\frac{3}{4}2^{\nu}}^{\frac{9}{8}2^{\nu+1}} \big( t^{-\sigma}\Vert  \Patm f \Vert_{p}\big)^q \frac{dt}{t}\Big]^{1/q}\lesssim \Big[\int\limits_{0}^{\infty} \big( t^{-\sigma}\Vert  \Patm f \Vert_{p}\big)^q \frac{dt}{t}\Big]^{1/q}$$
and this finishes the proof of (a).

Now we pass to the proof of (b). Note that
\begin{align*}
\Vert f\Vert_{\Bifi}^q=\sum\limits_{\eta\in\mathbb{Z}}\int\limits_{2^{\eta}}^{2^{\eta+1}}\big(t^{-\sigma}\Vert \Patm f\Vert_p\big)^q\frac{dt}{t}&\lesssim \sum\limits_{\eta\in\mathbb{Z}}\Big(2^{-\eta\sigma}\Big\Vert\sum\limits_{\nu\in\mathbb{Z}}\sum\limits_{Q\in\mathcal{D}_{\nu}} \vert s_Q\vert \sup\limits_{t\in \left(2^{\eta}, 2^{\eta+1}\right]} \vert \Patm  a_Q\vert \Big\Vert_p\Big)^q.
\end{align*}

Now we split the triple sum onto two sums over the sets $\{\nu:\ \nu>\eta\}$ and $\{\nu:\ \nu\leq\eta\}$ and denote the resulting sums by $I_1$ and $I_2$.

Fix $r<r_0$. Applying the assumptions we get $M>d/r-\sigma$ and $N>d(1/r-1)$ so $r>\frac{d}{d+N}$. Now, we apply Lemmas \ref{molest} and \ref{FrazJawLemma}. For $\nu >\eta$ we obtain
$$\sum\limits_{Q\in\mathcal{D}_{\nu}} \vert s_Q\vert \sup\limits_{t\in \left(2^{\eta}, 2^{\eta+1}\right]} \vert \Patm  a_Q(x)\vert\lesssim   2^{\nu\sigma}2^{(\eta-\nu)(m-N-d)}     \mathcal{M}_r\Big(\sum\limits_{Q\in\mathcal{D}_{\nu}}\vert s_Q\vert \vert Q\vert^{-1/p}\chi_Q  \Big)(x),$$
whereas for $\nu\leq\eta$ we get
$$\sum\limits_{Q\in\mathcal{D}_{\nu}} \vert s_Q\vert \sup\limits_{t\in \left(2^{\eta}, 2^{\eta+1}\right]} \vert \Patm  a_Q(x)\vert\lesssim 2^{\nu\sigma}2^{(-\eta+\nu)(M-d/r)}  \mathcal{M}_r\Big(\sum\limits_{Q\in\mathcal{D}_{\nu}}\vert s_Q\vert \vert Q\vert^{-1/p}\chi_Q  \Big)(x).$$
Hence, using the Minkowski inequality, the strong type of the maximal operator and finally the equality \eqref{FormulaFors_Q} we obtain 
\begin{align*}
I_1&\lesssim \sum\limits_{\eta\in\mathbb{Z}}\Big[\Big\Vert \sum\limits_{\nu>\eta} 2^{r(\eta-\nu)(m-N-d-\sigma)}  \Big[\mathcal{M}_r\Big(\sum\limits_{Q\in\mathcal{D}_{\nu}}\vert s_Q\vert \vert Q\vert^{-1/p}\chi_Q  \Big)\Big]^r    \Big\Vert_{p/r}^{1/r}\Big]^q\\
&\lesssim \sum\limits_{\eta\in\mathbb{Z}}\Big[ \sum\limits_{\nu>\eta} 2^{r(\eta-\nu)(m-N-d-\sigma)}  \Big\Vert\mathcal{M}_r\Big(\sum\limits_{Q\in\mathcal{D}_{\nu}}\vert s_Q\vert \vert Q\vert^{-1/p}\chi_Q  \Big)    \Big\Vert_{p}^r\Big]^{q/r}\\
&\lesssim\sum\limits_{\eta\in\mathbb{Z}}\Big( \sum\limits_{\nu>\eta} 2^{r(\eta-\nu)(m-N-d-\sigma)}\Big[\sum\limits_{Q\in\mathcal{D}_{\nu}}\vert s_Q\vert^p\Big]^{r/p}\Big)^{q/r},
\end{align*}
Now, Young's inequality and the assumption $m>N+d+\sigma $ yield
$$I_1\lesssim \sum\limits_{\nu\in\mathbb{Z}}\Big[\sum\limits_{Q\in\mathcal{D}_{\nu}}\vert s_Q\vert^p\Big]^{q/p}.$$
The proof for $I_2$ is similar: 
\begin{align*}
I_2&\lesssim \sum\limits_{\eta\in\mathbb{Z}}\Big( \sum\limits_{\nu\leq\eta} 2^{r(-\eta+\nu)(M-d/r+\sigma)}\Big[\sum\limits_{Q\in\mathcal{D}_{\nu}}\vert s_Q\vert^p \Big]^{r/p}\Big)^{q/r}.
\end{align*}
Young's inequality and the assumption $M>d/r-\sigma$ imply
$$I_2\lesssim \sum\limits_{\nu\in\mathbb{Z}}\Big[\sum\limits_{Q\in\mathcal{D}_{\nu}}\vert s_Q\vert^p \Big]^{q/p},$$
and this completes the proof.
\end{proof}
\begin{proof}[Proof of Theorem \ref{BesovEquiv}]
Suppose that the assumptions hold. Fix $N=\lceil d(1/r_0-1) \rceil$ and $M>\max\left(m_1,m_2,d/r_0-\sigma\right)$. If $f\in \dot{B}^{\sigma,\Lafi ,m_1}_{p,q}$, then Theorem \ref{moldecBesov} (a) implies that there exist a sequence $\{a_Q\}_{Q\in\mathcal{D}} $ of $(\Lafi, M, N, \sigma, p)$ molecules and a sequence of coefficients $\{s_Q\}_{Q\in\mathcal{D}}$  such that $f=\sum\limits_{Q\in\mathcal{D}}s_Qa_Q $ in $\Safid$ and
$$\Big[ \sum\limits_{\nu\in\mathbb{Z}}\Big(\sum\limits_{Q\in\mathcal{D}_{\nu}}\vert s_Q\vert^p \Big)^{q/p}\Big]^{1/q}\simeq \Vert f\Vert_{\dot{B}^{\sigma,\Lafi ,m_1}_{p,q} }. $$
Hence Theorem \ref{moldecBesov} (b) implies that $f\in \dot{B}^{\sigma,\Lafi ,m_2}_{p,q}$ and 
$$\Vert f\Vert_{\dot{B}^{\sigma,\Lafi ,m_2}_{p,q} }\lesssim \Big[ \sum\limits_{\nu\in\mathbb{Z}}\Big(\sum\limits_{Q\in\mathcal{D}_{\nu}}\vert s_Q\vert^p \Big)^{q/p} \Big]^{1/q}\simeq \Vert f\Vert_{\dot{B}^{\sigma,\Lafi ,m_1}_{p,q} }. $$
Similarly, if $f\in \dot{B}^{\sigma,\Lafi ,m_2}_{p,q}$, then $f\in \dot{B}^{\sigma,\Lafi ,m_1}_{p,q}$ and
$$\Vert f\Vert_{\dot{B}^{\sigma,\Lafi ,m_1}_{p,q} }\lesssim  \Vert f\Vert_{\dot{B}^{\sigma,\Lafi ,m_2}_{p,q} }. $$
\end{proof}

\begin{proof}[Proof of Theorem \ref{BesovEmbed}]
(i) Fix $m,N$ and $M$ satisfying the assumptions of Theorem \ref{moldecBesov}, where $r_0=\min(p,q_1,q_2,1)$. Let $f\in \dot{B}^{\sigma,\Lafi}_{p,q_1}$. Theorem  \ref{moldecBesov} (a) implies that there exist a sequence of $(\Lafi,M,N,\sigma,p)$ molecules $\{a_Q\}_{Q\in\mathcal{D}}$ and a sequence of coefficients  $\{s_Q\}_{Q\in\mathcal{D}}$ such that $f=\sum\limits_{Q\in\mathcal{D}}s_Qa_Q $ in $\Safid$ and
$$\Big[ \sum\limits_{\nu\in\mathbb{Z}}\Big(\sum\limits_{Q\in\mathcal{D}_{\nu}}\vert s_Q\vert^p \Big)^{q_1/p}\Big]^{1/q_1}\simeq \Vert f\Vert_{\dot{B}^{\sigma,\Lafi}_{p,q_1} }. $$
The inequality $q_1\leq q_2$ yields
$$\Big[ \sum\limits_{\nu\in\mathbb{Z}}\Big(\sum\limits_{Q\in\mathcal{D}_{\nu}}\vert s_Q\vert^p \Big)^{q_2/p}\Big]^{1/q_2} \lesssim \Big[ \sum\limits_{\nu\in\mathbb{Z}}\Big(\sum\limits_{Q\in\mathcal{D}_{\nu}}\vert s_Q\vert^p \Big)^{q_1/p}\Big]^{1/q_1}\simeq \Vert f\Vert_{\dot{B}^{\sigma,\Lafi}_{p,q_1} }. $$
Now, Theorem \ref{moldecBesov} (b) gives
$$\Vert f\Vert_{\dot{B}^{\sigma,\Lafi}_{p,q_2} }\lesssim\Vert f\Vert_{\dot{B}^{\sigma,\Lafi}_{p,q_1} }, $$
and the proof is completed.\\
$\ $\\
(ii) Note that a $(\Lafi, M, N, \sigma_1, p_1)$ molecule is a $(\Lafi, M, N, \sigma_2, p_2)$ molecule. Let $f\in \dot{B}^{\sigma_1,\Lafi}_{p_1,q}$. Theorem \ref{moldecBesov} (a) implies that there exist a sequence of $(\Lafi, M, N, \sigma_1, p_1)$ molecules $\{a_Q\}_{Q\in\mathcal{D}}$ and a sequence of coefficients $\{s_Q\}_{Q\in\mathcal{D}}$ such that $f=\sum\limits_{Q\in\mathcal{D}}s_Q a_Q$ in $\Safid$ and
$$\Big[\sum\limits_{\nu\in\mathbb{Z}}\Big(\sum\limits_{Q\in\mathcal{D}_{\nu}}\vert s_Q\vert ^{p_1} \Big)^{q/p_1}\Big]^{1/q}\lesssim \Vert f\Vert_{\dot{B}^{\sigma_1,\Lafi}_{p_1,q}}.$$
The inequality $p_2\geq p_1$ implies
$$\Big[\sum\limits_{\nu\in\mathbb{Z}}\Big(\sum\limits_{Q\in\mathcal{D}_{\nu}}\vert s_Q\vert ^{p_2}\Big)^{q/p_2}\Big]^{1/q}\lesssim \Big[\sum\limits_{\nu\in\mathbb{Z}}\Big(\sum\limits_{Q\in\mathcal{D}_{\nu}}\vert s_Q\vert ^{p_1}\Big)^{q/p_1}\Big]^{1/q}\lesssim \Vert f\Vert_{\dot{B}^{\sigma_1,\Lafi}_{p_1,q}}.$$
Hence Theorem \ref{moldecBesov} (b) implies
$$\Vert f\Vert_{\dot{B}^{\sigma_2,\Lafi}_{p_2,q}}\lesssim \Vert f\Vert_{\dot{B}^{\sigma_1,\Lafi}_{p_1,q}}, $$
and this finishes the proof.
\end{proof}

\section{Homogeneous Triebel-Lizorkin spaces associated to $\Lafi$}
In this section we introduce the Triebel-Lizor\-ki\-n spaces in the setting of Laguerre expansions of Hermite type, and state and prove the results, similar to these in Section 4. 
{\defin For $\sigma\in\mr$, $0<p<\infty,\ 0<q\leq\infty$, and $m\in\mathbb{N}$ such that $m>m_0$, we define the {\it homogeneous Triebel-Lizorkin spaces} $\TiLifi$ by
$$\TiLifi=\Big\{ f\in\Safid:\quad \Vert f\Vert_{\TiLifi}=\Big\Vert\Big(\int\limits_0^{\infty} \left\vert t^{-\sigma} \Patm f\right\vert^q\frac{dt}{t}\Big)^{1/q}\Big\Vert_p<\infty\Big\}.$$}\\

{\thm\label{MolDecTriebLiz}Let $\sigma\in\mr$, $0<p<\infty$ and $0<q\leq\infty$.
\begin{enumerate}[(a)]
\item For $M,N\in\mathbb{N}_+$ and $m>m_0$, if $f\in\TiLifi$, then there exist a sequence of $(\Lafi, M,N,\sigma,p)$ molecules $\{a_Q\}_{Q\in\mathcal{D}}$ and a sequence of coefficients $\{s_Q\}_{Q\in\mathcal{D}}$ such that $f=\sum\limits_{Q\in\mathcal{D}}s_Qa_Q $ in $\Safid$. Moreover,
\begin{equation}\label{SimCoeffTriebLiz}
\Big\Vert \Big[ \sum\limits_{Q\in\mathcal{D}}\left(\vert s_Q\vert \vert Q\vert^{-1/p}\chi_Q\right)^{q}\Big]^{1/q}\Big\Vert_p\simeq \Vert f\Vert_{\TiLifi} 
\end{equation}
\item Conversely, if  $M>\max(d/r_0-\sigma,m)$, $N>d(2/r_0-1)$, $m>\max(\sigma,0)+N+d$, and for a sequence of $(\Lafi, M,N,\sigma,p)$ molecules $\{a_Q\}_{Q\in\mathcal{D}}$ and a sequence of coefficients $\{s_Q\}_{Q\in\mathcal{D}}$ satisfying
$$\Big\Vert \Big[ \sum\limits_{Q\in\mathcal{D}}\big(\vert s_Q\vert \vert Q\vert^{-1/p}\chi_Q\big)^{q}\Big]^{1/q}\Big\Vert_p<\infty,$$
the series $\sum\limits_{Q\in\mathcal{D}}s_Qa_Q $ converges in $\Safid$, then its sum $f$ is in $\TiLifi$. Moreover,
$$\Vert f\Vert_{\TiLifi}\lesssim \Big\Vert \Big[ \sum\limits_{Q\in\mathcal{D}}\big(\vert s_Q\vert \vert Q\vert^{-1/p}\chi_Q\big)^{q}\Big]^{1/q}\Big\Vert_p.$$
\end{enumerate}}

\begin{proof} We begin with (a). Similarly as in the proof of Theorem \ref{moldecBesov} (a) we get $f=\sum\limits_{Q\in\mathcal{D}}s_Qa_Q $ in $\Safid$, where
$$s_Q=2^{-\nu\sigma}\vert Q\vert^{1/p}\sup\limits_{(y,t)\in Q\times (2^{\nu},2^{\nu+1}]} \left\vert \Patm f(y)\right\vert$$
and $a_Q=\left(\sqrt{\Lafi}\right)^M b_Q$, where
$$b_Q(x)=\frac{c_{m,M,N}}{s_Q}\int\limits_{2^{\nu}}^{2^{\nu+1}}t^M P^{\alpha}_{t,N} \left( \Patm f\cdot\chi_Q\right)(x)\frac{dt}{t},$$
and $a_Q$ are $(\Lafi, M,N,\sigma,p)$ molecules. Now, it suffices to prove \eqref{SimCoeffTriebLiz}. Firstly, we obtain
\begin{align*}
\Vert f\Vert_{\TiLifi}&=\Big\Vert\Big(\sum\limits_{\nu\in\mathbb{Z}}\int\limits_{2^{\nu}}^{2^{\nu+1}}\Big(\sum\limits_{Q\in\mathcal{D}_{\nu}} \vert t^{-\sigma} \Patm f\vert\chi_Q \Big)^q\frac{dt}{t}\Big)^{1/q}\Big\Vert_p\\
&\lesssim \Big\Vert \Big(\sum\limits_{\nu\in\mathbb{Z}}\sum\limits_{Q\in\mathcal{D}_{\nu}}\big( \vert Q\vert^{-1/p} \vert s_Q\vert\chi_Q\big)^q \Big)^{1/q}\Big\Vert_p.
\end{align*}
To prove the opposite inequality in \eqref{SimCoeffTriebLiz} note that as in Theorem \ref{moldecBesov} (a) we obtain 
$$\sum\limits_{Q\in\mathcal{D}_{\nu}}\vert Q\vert^{-1/p} \vert s_Q\vert \chi_Q(x)\lesssim \mathcal{M}_r \Big(\Big[ \int\limits_{\frac{3}{4}2^{\nu}}^{\frac{9}{8}2^{\nu+1}}\big( t^{-\sigma}\vert \Patm f\vert \big)^r \frac{dt}{t}\Big]^{1/r}\Big)(x),$$
where $0<r<r_0$. Hence, using the Fefferman-Stein vector valued maximal inequality and H\"{o}lder's inequality we get
\begin{align*}
\Big\Vert \Big[ \sum\limits_{\nu\in\mathbb{Z}}\sum\limits_{Q\in\mathcal{D}_{\nu}}\big( \vert Q\vert^{-1/p}s_Q\chi_Q\big)^q \Big]^{1/q}\Big\Vert_p &\lesssim \Big\Vert \Big[ \sum\limits_{\nu\in\mathbb{Z}} \Big( \int\limits_{\frac{3}{4}2^{\nu}}^{\frac{9}{8}2^{\nu+1}}\big( t^{-\sigma}\vert \Patm f\vert \big)^r \frac{dt}{t}   \Big)^{q/r}\Big]^{1/q}\Big\Vert_p\\
&\lesssim \Big\Vert \Big(\int\limits_0^{\infty} \big\vert t^{-\sigma} \Patm f\big\vert^q\frac{dt}{t}\Big)^{1/q}\Big\Vert_p
\end{align*}
and the last quantity equals $\Vert f\Vert_{\TiLifi}$. This completes the proof of (a).

For the proof of (b) we compute
\begin{align*}
\Vert f\Vert_{\TiLifi}&\lesssim \Big\Vert \Big(\sum\limits_{\eta\in\mathbb{Z}} \Big(\sum\limits_{\nu\in\mathbb{Z}}\sum\limits_{Q\in\mathcal{D}_{\nu}} 2^{-\eta\sigma} \vert s_Q\vert \sup\limits_{t\in[2^{\eta},2^{\eta+1}] }\vert\Patm a_Q\vert \Big)^q \Big)^{1/q} \Big\Vert_p\\
\end{align*}
Now we split the triple sum onto two sums over the sets $\{\nu:\ \nu>\eta\} $ and $\{\nu:\ \nu\leq\eta\} $ and denote the resulting sums by $I_1$ and $I_2$.

Similarly to the proof of Theorem \ref{moldecBesov} (b) for $\nu>\eta$ we get
$$\sum\limits_{Q\in\mathcal{D}_{\nu}}\vert s_Q\vert  \sup\limits_{t\in[2^{\eta},2^{\eta+1}] }\vert\Patm a_Q(x)\vert\lesssim 2^{\nu\sigma}2^{(\eta-\nu)(m-N-d)}\mathcal{M}_r \Big(\sum\limits_{Q\in\mathcal{D}_{\nu}}\vert s_Q\vert \vert Q\vert^{-1/p}\chi_Q  \Big)(x),$$
whereas for $\nu\leq\eta$
$$\sum\limits_{Q\in\mathcal{D}_{\nu}}\vert s_Q\vert \sup\limits_{t\in \left(2^{\eta}, 2^{\eta+1}\right]} \vert \Patm  a_Q(x)\vert\lesssim 2^{\nu\sigma}2^{(-\eta+\nu)(M-d/r)}  \mathcal{M}_r \Big(\sum\limits_{Q\in\mathcal{D}_{\nu}}\vert s_Q\vert \vert Q\vert^{-1/p}\chi_Q  \Big)(x).$$
Hence, applying Young's inequality we obtain
\begin{align*}
I_1&\lesssim \Big(\sum\limits_{\eta\in\mathbb{Z}}\Big[ \sum\limits_{\nu>\eta}2^{r(\eta-\nu)(m-N-d-\sigma)}\Big(\mathcal{M}_r \Big(\sum\limits_{Q\in\mathcal{D}_{\nu}}\vert s_Q\vert \vert Q\vert^{-1/p}\chi_Q  \Big)\Big)^r \Big]^{q/r}\Big)^{1/q}\\
&\lesssim \Big( \sum\limits_{\nu\in\mathbb{Z}}\Big[\mathcal{M}_r \Big(\sum\limits_{Q\in\mathcal{D}_{\nu}}\vert s_Q\vert \vert Q\vert^{-1/p}\chi_Q  \Big)\Big]^{q}\Big)^{1/q}.
\end{align*}
The Fefferman-Stein vector valued maximal inequality implies
$$\Vert I_1\Vert_p\lesssim \Big\Vert \Big[ \sum\limits_{\nu\in\mathbb{Z}}\Big(\sum\limits_{Q\in\mathcal{D}_{\nu}}\vert s_Q\vert \vert Q\vert^{-1/p}\chi_Q \Big)^{q}\Big]^{1/q}\Big\Vert_p= \Big\Vert \Big[ \sum\limits_{Q\in\mathcal{D}}\big(\vert s_Q\vert \vert Q\vert^{-1/p}\chi_Q\big)^{q}\Big]^{1/q}\Big\Vert_p.$$
In the same manner we deal with $I_2$:
\begin{align*}
I_2&\lesssim \Big(\sum\limits_{\eta\in\mathbb{Z}}\Big[ \sum\limits_{\nu\leq\eta}2^{(-\eta+\nu)(M-d/r+\sigma)}\mathcal{M}_r \Big(\sum\limits_{Q\in\mathcal{D}_{\nu}}\vert s_Q\vert \vert Q\vert^{-1/p}\chi_Q  \Big)\Big]^q \Big)^{1/q}\\
&\lesssim \Big( \sum\limits_{\nu\in\mathbb{Z}}\Big[\mathcal{M}_r \Big(\sum\limits_{Q\in\mathcal{D}_{\nu}}\vert s_Q\vert \vert Q\vert^{-1/p}\chi_Q  \Big)\Big]^{q}\Big)^{1/q},
\end{align*}
thus
$$\Vert I_2\Vert_p\lesssim \Big\Vert \Big[ \sum\limits_{Q\in\mathcal{D}}\left(\vert s_Q\vert \vert Q\vert^{-1/p}\chi_Q\right)^{q}\Big]^{1/q}\Big\Vert_p,$$
and this completes the proof.
\end{proof}

{\thm Let $\sigma\in\mr$. Assume that $0<p<\infty,\ 0<q\leq\infty$, $m_1,m_2\in\mathbb{N}$ and $m_1,m_2>m_0$. Then the spaces $\dot{F}^{\sigma,\Lafi ,m_1}_{p,q}$ and $\dot{F}^{\sigma,\Lafi ,m_2}_{p,q}$ coincide and their norms are equivalent.}\\

The proof is analogous to the proof of Theorem \ref{BesovEquiv} so we omit it.

\end{document}